# Chern numbers for singular varieties and elliptic homology

By Burt Totaro

A fundamental goal of algebraic geometry is to do for singular varieties whatever we can do for smooth ones. Intersection homology, for example, directly produces groups associated to any variety which have almost all the properties of the usual homology groups of a smooth variety. Minimal model theory suggests the possibility of working more indirectly by relating any singular variety to a variety which is smooth or nearly so.

Here we use ideas from minimal model theory to define some characteristic numbers for singular varieties, generalizing the Chern numbers of a smooth variety. This was suggested by Goresky and MacPherson as a next natural problem after the definition of intersection homology [11]. We find that only a subspace of the Chern numbers can be defined for singular varieties. A convenient way to describe this subspace is to say that a smooth variety has a fundamental class in complex bordism, whereas a singular variety can at most have a fundamental class in a weaker homology theory, elliptic homology. We use this idea to give an algebro-geometric definition of elliptic homology: "complex bordism modulo flops equals elliptic homology."

This paper was inspired by some questions asked by Jack Morava. The descriptions of elliptic homology given by Gerald Höhn [13] were also an important influence. Thanks to Dave Bayer, Mike Stillman, John Stembridge, Sheldon Katz, and Stein Stromme for their computer algebra programs Macaulay, SF, and Schubert, which helped in guessing the right answer.

## Contents

## 1. Statements

This paper presents two main results, which we will first state (using some new terminology) and then explain over the course of this section. First (Theorem 4.1), a rational linear combination of Chern numbers, viewed as an invariant of compact complex manifolds, is unchanged under "classical flops" if and only if it is a linear combination of the coefficients of the complex elliptic genus studied by Krichever and Höhn [19], [13]. This elliptic genus can be viewed as a power series associated to any compact complex manifold, the coefficients of the series being certain fixed linear combinations of the Chern numbers of the manifold. A more precise form of this result determines a geometrically meaningful version of complex elliptic homology over the ring $\mathbf{Z}[1/2]$ (Theorem 6.1, Remark 1). Second, we can ask when a given rational linear combination of Chern numbers, viewed as an invariant of smooth compact varieties, can be extended to an invariant of singular varieties, subject to a natural condition (compatibility with "IH-small resolutions"). This compatibility condition implies (as explained below) that the given linear combination of Chern numbers is invariant under classical flops; so every such linear combination of Chern numbers is a linear combination of the coefficients of the elliptic genus. We conjecture that, conversely, the elliptic genus can be defined for arbitrary singular varieties (compatibly with IH-small resolutions). The second main result of this paper is that at least a certain weaker invariant, which Höhn called the twisted $\chi_y$ genus, can be defined for arbitrary singular varieties (compatibly with IH-small resolutions): see Theorems 8.1 and 8.2.

From now on, we call any homogeneous rational polynomial in variables $c_1, \ldots, c_n$ of degree $n$ a Chern number (for $n$-folds). Here the variable $c_i$ is given degree $i$. A Chern number in this sense determines a function from compact complex manifolds of complex dimension $n$ to the rational numbers: replace $c_1, \ldots, c_n$ by the Chern classes of the tangent bundle and integrate the resulting top-degree cohomology class. Moreover, a polynomial in $c_1, \ldots, c_n$ is uniquely determined by its values on compact complex manifolds, or even just on smooth complex projective varieties. (The reason is that Chern numbers can be identified with linear functions $MU_{2n} \otimes \mathbf{Q} \to \mathbf{Q}$, where $MU_{2n}$ is the bordism group of weakly complex $2n$-manifolds (see Section 2), and the group $MU_{2n}$ is generated by smooth complex projective $n$-folds [23].) For example, the Euler characteristic, the signature, and the Todd genus, for smooth projective varieties of a given dimension, are Chern numbers, thanks to Hirzebruch [12].

Let us define an IH-small resolution of a singular variety $Y$ to be a resolution of singularities $f : X \to Y$ such that for every $i \geq 1$, the set of points $y \in Y$ such that $\dim(f^{-1}(y)) = i$ has codimension greater than $2i$ in $Y$. The interest of such a resolution is that there is a natural identification of the intersection homology of $Y$ with the ordinary homology of $X$, by [10, p. 121].

On the other hand, the intersection homology of $Y$ is a direct summand in the homology of any resolution of singularities of $Y$; so IH-small resolutions, when they exist, are the "smallest possible" resolutions of a given variety. In fact, IH-small resolutions turn out to be relative minimal models in the precise sense of Mori's program (see Section 8), and this is crucial to our approach. We remark that most singular varieties have no IH-small resolution.

The problem we are considering was formulated by Goresky and MacPherson [11, Problem 10]: Which Chern numbers $\alpha$ for $n$-folds can be defined for all singular compact complex $n$-folds $Y$ in such a way that, whenever $f : X \to Y$ is an IH-small resolution, we have $\alpha(Y) = \alpha(X)$? Notice that $\alpha(X)$ is already defined since $X$ is smooth. Also, the question makes sense either for $X$, $Y$ compact complex spaces or for $X$, $Y$ projective varieties; we usually assume $X$ and $Y$ projective in this paper, although we will say when something works more generally. (We expect the same answer in both situations.)

When we say that a given Chern number "can be defined" for singular varieties, we always mean that it can be defined compatibly with IH-small resolutions in the above sense. The motivation for this problem is that, as Goresky and MacPherson observed, intersection homology provides definitions of the Euler characteristic, the signature, and the Todd genus for all singular varieties, and the resulting definitions are compatible with IH-small resolutions in this sense.

To define more general Chern numbers, one might hope to lift the homology Chern classes of a singular variety to intersection homology and then multiply some of them. Indeed, as Goresky and MacPherson conjectured, all algebraic cycles, and in particular the homology Chern classes, lift rationally to intersection homology, by Barthel, Brasselet, Fieseler, Gabber, and Kaup [3]; but those lifts are not unique (see the comments on this problem in [3, p. 158]), and we will not use that approach.

Clearly, if a given Chern number can be defined for singular varieties in the above sense, then it must take the same value on any two IH-small resolutions of a given singular variety. Thus, in order to give an upper bound for the rational vector space of degree-$n$ polynomials in $c_1, \ldots, c_n$ which can be defined for singular varieties, we need an explicit collection of singular varieties $Y$ which have two different IH-small resolutions. We will use those singular varieties $Y$ of dimension $n \geq 3$ which are Zariski locally isomorphic, near each point of their singular set $Z$, to the product of a 3-fold node with a smooth $(n-3)$-fold (see Section 4 for more details). Such a variety has two different IH-small resolutions $X_1$ and $X_2$, and in this situation we say that $X_1$ and $X_2$ are related by a *classical flop*. (The manifold $X_2$ is obtained from $X_1$ by cutting out a $\mathbf{P}^1$-bundle over $Z$ and replacing it by a possibly different $\mathbf{P}^1$-bundle over $Z$.)

Now at last we can understand the first main result of this paper, Theorem 4.1: the space of Chern numbers which do not change under classical flops is spanned by the coefficients of the complex elliptic genus. This genus is a homomorphism of graded rings from the complex bordism ring

$$MU_* \otimes \mathbf{Q} = \mathbf{Q}[\mathbf{CP}^1, \mathbf{CP}^2, \ldots]$$

onto

$$\mathrm{Ell}_* := \mathbf{Q}[x_1, x_2, x_3, x_4],$$

where we put the bordism group $MU_{2n}$ in degree $n$ for all $n$, and $x_n$ in degree $n$ for $1 \leq n \leq 4$ [19], [13]. (The more famous elliptic genus defined by Landweber, Stong, Ochanine, and Witten [20] is the homomorphism from $MU_* \otimes \mathbf{Q}$ to $\mathbf{Q}[x_2, x_4]$ defined by setting $x_1$ and $x_3$ to 0 in the complex elliptic genus; the resulting invariant is then defined for oriented manifolds, not just complex manifolds.)

A reformulation of our first main result is that "complex bordism modulo classical flops equals elliptic homology," at least rationally. In fact, a suitable formulation of this statement is true over $\mathbf{Z}[1/2]$, by a more refined version, Theorem 6.1, of our theorem.

In the complex elliptic genus, the coefficient of each monomial in $x_1, \ldots, x_4$ of degree $n$ is a certain Chern number of degree $n$. (We repeat that in this paper "Chern number" means a homogeneous rational polynomial in variables $c_1, c_2, \ldots$ of degrees $1, 2, \ldots$. The coefficients of the complex elliptic genus in the sense just stated are invariant under classical flops, but the individual Chern monomials which occur in these coefficients are generally not invariant under classical flops.) It follows from our first main result that the space of Chern numbers which can be defined for all singular varieties is *at most* equal to those given by the complex elliptic genus. The main problem left open by this paper is to show that this is an equality, that is, to define the elliptic genus of an arbitrary singular variety.

One partial result in that direction follows immediately from earlier work. Thanks to Morihiko Saito's Hodge structure on intersection homology [32], [33], it was already known how to define Hodge numbers for a singular variety in a way which is compatible with IH-small resolutions. This immediately defines a few Chern numbers for singular varieties, namely those corresponding to the Hirzebruch $\chi_y$ genus [12]. We can view the Hirzebruch $\chi_y$ genus as a surjective homomorphism of graded rings,

$$\chi_y : MU_* \otimes \mathbf{Q} \to \mathbf{Q}[x_1, x_2].$$

The $\chi_y$ genus includes the Euler characteristic, the signature, and the Todd genus as special cases.

The second main result of this paper is that the space of Chern numbers which take the same value on any two IH-small resolutions of a given singular variety is larger than just the $\chi_y$ genus (Theorem 8.1). Namely, Höhn defined a genus, the twisted $\chi_y$ genus, which is intermediate between the $\chi_y$ genus and the complex elliptic genus [13, p. 65]. For a smooth variety $X$, the twisted $\chi_y$ genus is defined as the set of holomorphic Euler characteristics of the bundles $\Omega^i_X \otimes K_X^{\otimes j}$, for $0 \leq i \leq n$ and $j \in \mathbf{Z}$. As he remarks, all the "classical" genera such as the signature, Todd genus, Euler characteristic, and the $\widehat{A}$ genus factor through the twisted $\chi_y$ genus. We show that the twisted $\chi_y$ genus takes the same value on any two IH-small resolutions of a given singular variety. As Höhn states, the twisted $\chi_y$ genus can be viewed as a surjective homomorphism

$$\chi_{yz} : MU_* \otimes \mathbf{Q} \to \mathbf{Q}[x_1, x_2, x_3, x_4]/(\Delta(x_2, x_3, x_4)),$$

where $\Delta$ is the expression for the discriminant cusp form $\Delta$ in the theory of modular forms as a polynomial in certain explicit Jacobi forms: $x_2$ is 24 times the Weierstrass $\mathfrak{p}$-function, $x_3$ is the derivative of the Weierstrass $\mathfrak{p}$-function, and $x_4 = 6\mathfrak{p}^2 - g_2/2$ where $g_2$ is the Eisenstein series of weight 4. See Section 9 for details.

In Section 8, we define an explicit extension of the twisted $\chi_y$ genus from smooth varieties to all singular varieties which possess a relative canonical model. The extension is compatible with IH-small resolutions when they exist. The minimal model conjecture would imply that every singular variety has a relative canonical model [17].

For example, for $n \leq 4$ the twisted $\chi_y$ genus includes all Chern numbers for $n$-folds. So we know that all Chern numbers for $n$-folds with $n \leq 4$ take the same value on any two IH-small resolutions of a singular variety. For all $n \leq 11$, the twisted $\chi_y$ genus and the complex elliptic genus are equal, so in these dimensions we know exactly which Chern numbers can be defined for singular varieties compatibly with IH-small resolutions. For example, the space of Chern numbers which can be defined for singular 5-folds happens to be spanned by certain Chern monomials, namely all of the Chern monomials except $c_3c_2$ (that is, $c_5$, $c_4c_1$, $c_3c_1^2$, $c_2^2c_1$, $c_2c_1^3$, and $c_1^5$). In fact, Goresky and MacPherson gave an example in [11, p. 222], of a 5-dimensional Schubert variety with two different IH-small resolutions, and I was led to the results of this paper by computing that these two resolutions have the same Chern numbers $c_5$, $c_4c_1$, $c_3c_1^2$, $c_2^2c_1$, $c_2c_1^3$, and $c_1^5$ as each other, but different $c_3c_2$'s.

We make some additional remarks. Any polynomial in the Chern classes of degree $n$ gives not only a Chern number for all smooth $n$-folds, but also a homology class (in $H_{2k}(X, \mathbf{Q})$) for all smooth $(n+k)$-folds, $k \geq 0$, and we can ask which of these classes can be defined as homology classes on all singular varieties in a way compatible with IH-small resolutions. It appears that the answer should be exactly the same as for Chern numbers. In particular, we can define homology classes corresponding to the twisted $\chi_y$ genus. If we could do this for the complex elliptic genus, we would have a natural fundamental class in rational elliptic homology $\mathrm{Ell}_{2n}(Y)$, for any singular compact $n$-fold $Y$. (Since we are working rationally, we can define elliptic homology here as a quotient of complex bordism: $\mathrm{Ell}_*(Y) := MU_*(Y) \otimes_{MU_*} \mathrm{Ell}_*$.) Notice that we should probably not expect to have well-defined characteristic classes for singular varieties in intersection homology as opposed to ordinary homology, since even the simplest Chern polynomials, namely the Chern classes, can be different in two different IH-small resolutions of the same singular variety, as Verdier found [5].

Also, there is a natural integral version of the question we have been considering rationally. One could try to compute the quotient ring of the complex bordism ring $MU_*$ by flops, but this ring seems not so natural integrally; for example, it is not finitely generated, although after tensoring with $\mathbf{Q}$ it becomes the polynomial ring $\mathbf{Q}[x_1, x_2, x_3, x_4]$. The natural integral question seems to be to compute the quotient ring of the $SU$-bordism ring $MSU_*$ by "$SU$-flops." Away from the prime 2, this works beautifully: the quotient of $MSU_* \otimes \mathbf{Z}[1/2]$ by $SU$-flops is equal to the image of $MSU_* \otimes \mathbf{Z}[1/2]$ under the complex elliptic genus, and this image is a polynomial ring $\mathbf{Z}[1/2][x_2, x_3, x_4]$ (Theorem 6.1). This analysis also determines a natural version of complex elliptic homology over $\mathbf{Z}[1/2]$, which is defined for the first time here. All this suggests the possibility of defining some version of elliptic homology as bordism with respect to some natural class of singular spaces which would include all Gorenstein complex varieties.

## 2. Weakly complex manifolds

We make some elementary remarks about weakly complex manifolds, the objects used to define the complex bordism ring $MU_*$, for use in Section 4. A weakly complex manifold is defined to be a real (smooth) manifold with a complex structure on the stable tangent bundle. More explicitly, given a real manifold $X$, its tangent bundle determines a homotopy class of maps $X \to BO$, and a weakly complex structure on $X$ is a homotopy class of lifts:

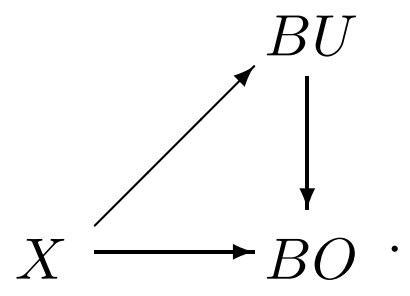

In particular, a complex manifold is a weakly complex manifold in a natural way. Likewise, an $SU$-manifold is defined to be a real manifold together with a homotopy class of lifts of the tangent bundle to $BSU$.

Here we say that two lifts $X \to BU$ of the tangent bundle $X \to BO$ are homotopic if they are homotopic through lifts; it is not enough for them to be homotopic just as maps $X \to BU$. In particular, there are two different weakly complex structures on a point, corresponding to $\pi_0(O/U) = \mathbf{Z}/2$, one coming from the complex structure on a point and the other not. In the bordism ring $MU_*$, the ring of closed weakly complex manifolds modulo boundaries of compact weakly complex manifolds [23], the first of these weakly complex structures on a point represents 1 and the second represents $-1$.

As a result, for any weakly complex manifold $X$, we can form a new weakly complex structure on the same real manifold which we call the negative weakly complex manifold, or $-X$, by taking the product of $X$ with the nontrivial weakly complex structure on a point. If the manifold $X$ is closed, then $-X$ is indeed the negative of $X$ in the bordism ring $MU_*$. Even if $X$ is a complex manifold, $-X$ is in general only a weakly complex manifold; the Chern classes of $-X$ are the same as those of $X$, but $-X$ has the opposite orientation.

The following lemma, which we will use in the proof of Theorem 4.1, is probably well known. It works for any version of bordism (unoriented, oriented, and so on).

Lemma 2.1. *Let $A$, $B$, $C$ be compact weakly complex manifolds. Suppose there are given diffeomorphisms of the boundaries of $A$, $B$, and $C$ to the same weakly complex manifold $M$. Then*

$$A \cup_M -B + B \cup_M -C + C \cup_M -A = 0$$

*in the bordism group $MU_*$.*

*Proof.* Let $H$ denote a hexagon in $\mathbf{R}^2$ (thus $H$ is homeomorphic to the disk). Let $W$ be the union of $A \times [0,1]$, $B \times [0,1]$, $C \times [0,1]$, and $M \times H$ modulo

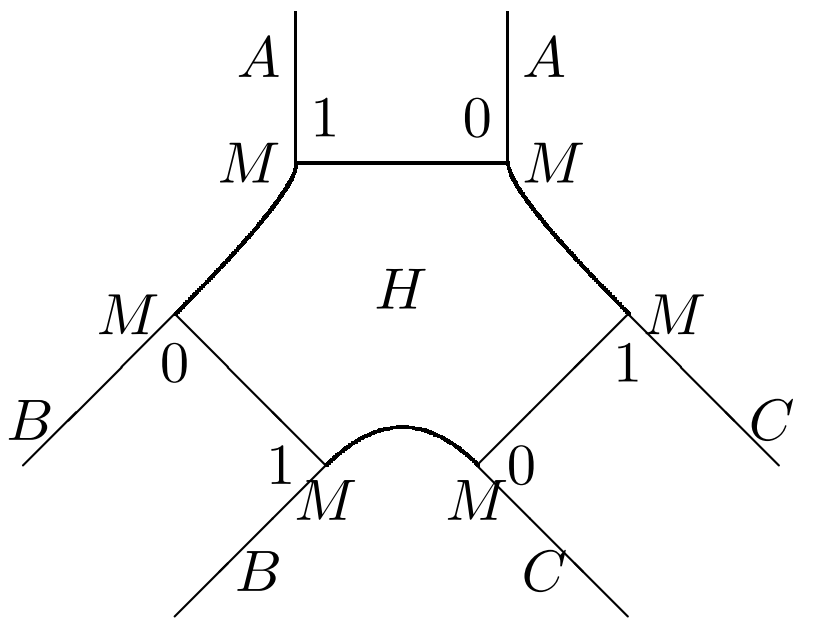

the identifications pictured here.

Then $W$ is a compact weakly complex manifold whose boundary is the disjoint union of the closed weakly complex manifolds $A \cup_M -B$, $B \cup_M -C$, and $C \cup_M -A$. □

## 3. The complex elliptic genus

We define the complex elliptic genus

$$MU_* \to \mathbf{Q}[x_1, x_2, x_3, x_4]$$

using one of the approaches in Höhn's thesis (Section 2.5) [13]. As it happens, the definition will not be used explicitly in most of this paper; the important thing is the *rigidity* property of this genus, which we will state in Section 4.

For a complex vector bundle $E$, define

$$\Lambda_t(E) = \oplus \Lambda^k E \cdot t^k$$

and

$$S_t(E) = \oplus S^k E \cdot t^k$$

as the power series in $t$ whose coefficients are the exterior or symmetric powers of $E$. Clearly $S_t(E+F) = S_t(E)S_t(F)$ and likewise for $\Lambda_t$. Also, these operations extend to virtual bundles, with $S_t(-E) = \Lambda_{-t}(E)$.

For complex numbers $\tau$ in the upper half-plane and $z \in \mathbf{C}$, the Weierstrass sigma function is defined by

$$\sigma(\tau, z) = z \prod_{\substack{\omega \in \mathbf{Z}+\mathbf{Z}\tau \\ \omega \neq 0}} (1 - z/\omega) e^{z/\omega + (z/\omega)^2/2}.$$

This is an entire function of $z$ with zero set equal to the lattice $\mathbf{Z} + \mathbf{Z}\tau \subset \mathbf{C}$. If we modify the sigma function slightly by defining

$$\Phi(\tau, z) = 2\pi i e^{\eta z^2/2 - \pi i z} \sigma(\tau, z),$$

we get a function which is periodic under $\tau \mapsto \tau + 1$ as well as $z \mapsto z + 1$ for a unique function $\eta(\tau)$. So $\Phi(\tau, z)$ admits an expansion in $q := e^{2\pi i \tau}$ and $y := e^{2\pi i z}$ [21, p. 247]:

$$\Phi(\tau, z) = \prod_{m \geq 1} \frac{(1 - y^{-1} q^{m-1})(1 - y q^m)}{(1 - q^m)^2}.$$

Our conventions are slightly different from Höhn's here: the function he calls $\Phi(\tau, z)$ vanishes for $z$ in the lattice $2\pi i(\mathbf{Z}+\mathbf{Z}\tau)$, rather than the more traditional $\mathbf{Z} + \mathbf{Z}\tau$ as here [13, p. 59]. Also, he writes $y = -e^z$ rather than our $y = e^{2\pi i z}$, so the definition below of the complex elliptic genus of a complex manifold differs from Höhn's by replacement of $y$ by $-y$.

We define the complex elliptic genus as the ring homomorphism

$$MU_* \to \mathbf{Q}((y))[[q,k]]$$

associated, in the way recalled below, to the following Hirzebruch characteristic power series $Q(x) \in \mathbf{Q}((y))[[q,k,x]]$ [13, p. 59]:

$$\begin{aligned} Q(x) &= e^{kx}\frac{x\Phi(\tau, x/2\pi i - z)}{\Phi(\tau, x/2\pi i)\Phi(\tau,-z)} \\ &= \Phi(q,1/y)^{-1}\frac{x}{1-e^{-x}}e^{kx}\prod_{m\geq 1}\frac{(1-yq^{m-1}e^{-x})(1-y^{-1}q^m e^x)}{(1-q^m e^x)(1-q^m e^{-x})}. \end{aligned}$$

As above, we write $q = e^{2\pi i\tau}$ and $y = e^{2\pi i z}$. Since the factor $\Phi(q,1/y)^{-1}$ does not involve $x$, it could be omitted without changing the important properties of the genus.

In Hirzebruch's general correspondence between power series and genera, we think of the variable $x$ as the first Chern class of a line bundle; then a series $Q(x) \in R[[x]]$, for a $\mathbf{Q}$-algebra $R$, gives a characteristic class for line bundles with values in $H^*(\cdot, R)$. This extends uniquely to an exponential characteristic class $\varphi(E)$ for arbitrary vector bundles $E$, where "exponential" means that $\varphi(E+F) = \varphi(E)\varphi(F)$. Then we get a genus $\varphi : MU_* \to R$ by associating, to a compact complex manifold $X$, the element $\varphi(X) := \int_X \varphi(TX)$ of $R$.

Let us work out what this means for the above series $Q(x)$. If $X$ is a compact complex $n$-manifold, let $x_1, \ldots, x_n$ denote the Chern roots of the tangent bundle $TX$. (These are formal variables whose symmetric functions are the Chern classes of $TX$.) Then

$$\begin{aligned} \varphi(X) &= \int_X Q(x_1)\cdots Q(x_n) \\ &= \int \mathrm{td}(TX)\mathrm{ch}(K_X^{\otimes -k}\otimes \prod_{m\geq 1}(\Lambda_{-y^{-1}q^m}T\otimes\Lambda_{-yq^{m-1}}T^*\otimes S_{q^m}T\otimes S_{q^m}T^*)), \end{aligned}$$

where $T = TX - n$. Here the Todd genus $\mathrm{td}(TX)$ comes from the factor

$$\prod_{i=1}^n x_i/(1-e^{-x_i})$$

in $Q(x_1)\cdots Q(x_n)$, and the factor $\Phi(q,1/y)^{-1}$ in $Q(x)$ merely has the effect of replacing $TX$ by the rank-0 virtual bundle $T$ in the above expression. By the Hirzebruch-Riemann-Roch theorem, we deduce that the complex elliptic genus of a compact complex manifold $X$ is given by the holomorphic Euler characteristic

$$\varphi(X) = \chi(X, K_X^{\otimes -k} \otimes \prod_{m \geq 1} (\Lambda_{-y^{-1}q^m} T \otimes \Lambda_{-yq^{m-1}} T^* \otimes S_{q^m} T \otimes S_{q^m} T^*))$$

in the ring $\mathbf{Q}((y))[[q,k]]$. Because of our slightly different conventions, this series differs by replacing $y$ by $-y$ from Höhn's definition of the genus.

The complex elliptic genus has better properties if $X$ is an $SU$-manifold, that is, if the canonical line bundle $K_X = \Lambda^n(T^*X)$ is trivial. For $X$ an $SU$-manifold, the series $\varphi(X)$ clearly lies in $\mathbf{Z}((y))[[q]]$. Moreover, for $X$ an $SU$-manifold of complex dimension $n$, $\varphi(X)$ is in fact a Jacobi form of weight $n$, by Höhn [13]. Jacobi forms are generalizations of modular forms defined by Eichler and Zagier [8], although we use a slight variant of their definition as we will explain later in this section. Just as modular forms (of level 1) are exactly sections of powers of a certain line bundle $\psi_1$ on the compactified moduli stack $\overline{M}_{1,1}$ of elliptic curves, we define a Jacobi form of weight $n$ to be a section of $n$ times a certain line bundle $\psi_1$ on the universal elliptic curve $\overline{M}_{1,2}$. In fact, this approach gives a definition of the ring of Jacobi forms with coefficients in any given commutative ring $R$. Here $\overline{M}_{g,r}$ is the Knudsen-Deligne-Mumford moduli stack of $r$-pointed stable curves of genus $g$, which comes with $r$ line bundles $\psi_1, \ldots, \psi_r$ representing the cotangent line of the curve at the $r$ given points [15]. The line bundle $\psi_1$ on $\overline{M}_{1,2}$ is not the pullback of the line bundle $\psi_1$ on $\overline{M}_{1,1}$ by the projection $\pi : \overline{M}_{1,2} \to \overline{M}_{1,1}$, forgetting the second point; instead, we have $\psi_1 = \pi^*\psi_1 + D_{0,1}$, where $D_{0,1}$ is the divisor on $\overline{M}_{1,2}$ where the two points lie on the same genus-0 component, or equivalently it is the zero section of the universal elliptic curve $\pi : \overline{M}_{1,2} \to \overline{M}_{1,1}$ [15]. The Tate elliptic curve $E$ [21, pp. 197–198] is a particular family of 1-pointed stable curves of genus 1 over $\mathbf{Z}[[q]]$, so there is a map $E \to \overline{M}_{1,2}$. By pulling back, a Jacobi form over $R$ restricts to a power series in $R((y))[[q]]$, when we use the natural coordinate $y$ on $E$. This restriction map is clearly injective. It is in this sense that the complex elliptic genus of an $SU$-manifold in $\mathbf{Z}((y))[[q]]$ is a Jacobi form over the integers.

The ring of (level 1) modular forms over any $\mathbf{Z}[1/6]$-algebra $R$ is the polynomial ring $R[g_2, g_3]$, where $g_i$ is the Eisenstein series of weight $2i$, $i \geq 2$. The ring is more complicated in characteristics 2 and 3 [6]. By the methods of Deligne's paper, one can also compute the ring of Jacobi forms, in the above sense. For any $\mathbf{Z}[1/6]$-algebra $R$, one gets the polynomial ring $R[x, y, g_2]$, where $x$ is the Weierstrass $\mathfrak{p}$-function (which has weight 2), $y$ is its derivative (of weight 3), and $g_2$ is the Eisenstein series of weight 4: this contains the ring $R[g_2, g_3]$ of modular forms, thanks to the Weierstrass equation [6, p. 59]:

$$y^2 = 4x^3 - g_2 x - g_3.$$

The series expansions of these Jacobi forms are as follows (Lang [21]), except that Lang puts an extra factor of $(2\pi i)^n$ in his definition of each Jacobi form of weight $n$:

$$\begin{aligned}
x(\tau,z) &= \frac{1}{12}+\frac{y}{(1-y)^2}-2\sum_{m,n\geq 1} nq^{mn}+\sum_{m,n\geq 1} nq^{mn}(y^n+y^{-n}),\\
y(\tau,z) &= \sum_{m\geq 0}\frac{q^m y(1+q^m y)}{(1-q^m y)^3}-\sum_{m\geq 1}\frac{(q^m/y)(1+q^m/y)}{(1-q^m/y)^3},\\
g_2(\tau) &= \frac{1}{12}\left[1+240\sum_{m\geq 1}\frac{m^3q^m}{1-q^m}\right],\\
g_3(\tau) &= \frac{1}{6^3}\left[-1+504\sum_{m\geq 1}\frac{m^5q^m}{1-q^m}\right].
\end{aligned}$$

The Jacobi forms $x$ and $g_2$ are only defined over $\mathbf{Z}[1/6]$ (equivalently, the coefficients of the corresponding power series have denominators), but if we let $x_2 = 24x$, $x_3 = y$, and $x_4 = 6x^2 - g_2/2$, then these Jacobi forms are defined over the integers. Using the methods of Deligne's paper [6], one finds that over any $\mathbf{Z}[1/2]$-algebra $R$, the ring of Jacobi forms is the polynomial ring $R[x_2, x_3, x_4]$.

Moreover, Höhn showed that the Jacobi forms $x_2, x_3, x_4$ arise as the elliptic genera of certain explicit $SU$-manifolds, of complex dimensions 2, 3, 4: the K3 surface, the almost complex 6-sphere, and a nonstandard weakly complex structure on a quadric 4-fold [13, pp. 24–25]. It follows that the ring homomorphism

$$MSU_* \to (\text{Jacobi forms over } \mathbf{Z})$$

becomes surjective after it is tensored with $\mathbf{Z}[1/2]$. Höhn pointed out that it is not surjective integrally: the Jacobi form $x_2$, 24 times the Weierstrass $\mathfrak{p}$-function, is the elliptic genus of the K3 surface, which generates $MSU_{2n}$ for $n = 2$: but $12\mathfrak{p}$, not only $24\mathfrak{p}$, is a Jacobi form with integer coefficients.

For clarity, let us explain the relation of this definition of Jacobi forms to Eichler and Zagier's slightly different analytic notion of Jacobi forms. In their terminology, the series $\varphi(X)$ associated to an $SU$-manifold of complex dimension $n$ is a meromorphic Jacobi form of weight $n$ and index 0 which is holomorphic outside the lattice $z \in \mathbf{Z}+\mathbf{Z}\tau$. We now give an analytic definition of Jacobi forms in our sense and explain the precise relation to Eichler and Zagier's definitions below. Namely, we call a power series $\varphi(q, y) \in \mathbf{C}((y))[[q]]$ a Jacobi form of weight $n$ if it has the following properties. It converges for $y$ sufficiently close to 0 and not equal to 0, and $q$ sufficiently close to 0 depending on $y$. It extends to a meromorphic function $\varphi(q, y)$ on $D \times (\mathbf{C} - 0)$, where $D$ is the unit disk, which is holomorphic outside the divisors $y = q^m$, $m \in \mathbf{Z}$. Changing variables by $q = e^{2\pi i\tau}$ and $y = e^{2\pi iz}$, we get a meromorphic function $\varphi(\tau, z)$ on $H \times \mathbf{C}$, where $H$ is the upper half-plane, which is holomorphic except for $z$ in the lattice $\mathbf{Z} + \mathbf{Z}\tau$. It satisfies

$$\varphi(\tau, z + \omega) = \varphi(\tau, z) \text{ for all } \omega \in \mathbf{Z} + \mathbf{Z}\tau;$$

that is, $\varphi(\tau, \cdot)$ is an elliptic function with respect to the lattice $\mathbf{Z} + \mathbf{Z}\tau \subset \mathbf{C}$. And it satisfies

$$\varphi(\frac{a\tau + b}{c\tau + d}, \frac{z}{c\tau + d})(c\tau + d)^{-n} = \varphi(\tau, z)$$

for all

$$\begin{pmatrix} a & b \\ c & d \end{pmatrix} \in SL(2, \mathbf{Z}).$$

We can check that the ring of Jacobi forms in the above sense coincides with the algebraically defined ring

$$\oplus_n H^0((\overline{M}_{1,2})_{\mathbf{C}}, n\psi_1) = \mathbf{C}[\mathfrak{p}, \mathfrak{p}', g_2]$$

by comparing with the results of Eichler and Zagier. A Jacobi form of weight $k$ in the above sense fails to be a "weak Jacobi form of weight $k$ and index 0" in their sense [8, p. 104] only because of its possible poles for $y = q^m$, $m \in \mathbf{Z}$. But there is a particular weak Jacobi form $\widetilde{\varphi}_{-2,1}$ of weight $-2$ and index 1 which vanishes to order 2 for $y = q^m$, $m \in \mathbf{Z}$, and nowhere else on $D \times (\mathbf{C} - 0)$ [8, p. 108]. Up to simple factors, $\widetilde{\varphi}_{-2,1}$ is the square of the Weierstrass sigma function; more precisely, in terms of the normalization $\Phi(q, y)$ of the sigma function used earlier in this section, we have

$$\widetilde{\varphi}_{-2,1}(\tau, z) = y\Phi(q, y)^2.$$

It follows that the ring of Jacobi forms in the above sense is the ring $\widetilde{J}_{*,*}[\widetilde{\varphi}_{-2,1}^{-1}]_{*,0}$, where $\widetilde{J}_{*,*}$ is the ring of weak Jacobi forms, bigraded by weight and index, and we only look at the localized ring in index 0. Now Eichler and Zagier's computation of the ring $\widetilde{J}_{*,*}$ [8, pp. 111–112] shows that the above localized ring is the polynomial ring $\mathbf{C}[\mathfrak{p}, \mathfrak{p}', g_2]$ (where $g_2$ is called $E_4$ in Eichler and Zagier, up to a constant factor). We deduce that the ring of Jacobi forms as defined in the previous paragraph is exactly the polynomial ring $\mathbf{C}[\mathfrak{p}, \mathfrak{p}', g_2]$.

The following fact, at least over $\mathbf{C}$, is implicit in Eichler and Zagier's calculation. We give a direct proof since we will need this later.

LEMMA 3.1. *Consider the universal elliptic curve $\overline{M}_{1,2}$ as a smooth stack* (*or orbifold*) *over a field k. Then any rational section of the line bundle $a\psi_1$, $a \in \mathbf{Z}$, over $\overline{M}_{1,2}$ which is regular outside the zero section $D_{0,1}$ is regular everywhere.* (*So it is a Jacobi form over k.*)

*Proof.* It suffices to show that for all integers $a$ and all $b \geq 1$, every (regular) section of the line bundle $a\psi_1 + bD_{0,1}$ over the surface $\overline{M}_{1,2}$ vanishes on the curve $D_{0,1}$. This follows if we can show that this line bundle has negative degree on $D_{0,1}$, which it does because $\psi_1 \cdot D_{0,1} = 0$ and $D_{0,1}^2 < 0$. The first statement follows from the more precise fact that the line bundle $\psi_1$ has trivial restriction to $D_{0,1}$: if we think of $\psi_1$ as the line bundle of 1-forms on a stable curve of genus 1 which have at most a pole at the origin, then a trivialization along the zero section is given by the residue. Since $D_{0,1}$ maps isomorphically to $\overline{M}_{1,1}$ under the projection $\pi : \overline{M}_{1,2} \to \overline{M}_{1,1}$, and the line bundle called $\psi_1$ on $\overline{M}_{1,1}$ has positive degree (in fact, degree 1/24; see below), the pullback line bundle $\pi^*\psi_1$ on $\overline{M}_{1,2}$ has positive degree on $D_{0,1}$. Since $\psi_1 = \pi^*\psi_1 + D_{0,1}$ on $\overline{M}_{1,2}$, it follows that $D_{0,1}^2 < 0$. Essentially the same calculation, together with some background on stacks, can be found in Mumford [27, p. 326], where he shows that

$$\int_{\overline{M}_{1,1}} \psi_1 = \int_{\overline{M}_{1,2}} \psi_1^2 = 1/24. \qquad \square$$

We return to topology by describing the image of the complex elliptic genus on arbitrary complex manifolds rather than on $SU$-manifolds. Namely, the bordism ring $MU_* \otimes \mathbf{Q}$ is a polynomial ring over $MSU_* \otimes \mathbf{Q}$ generated by $\mathbf{CP}^1$, by Novikov [29]. And it is easy to check that the series $\varphi(\mathbf{CP}^1) \in \mathbf{Q}((y))[[q,k]]$ is algebraically independent of the image of $MSU_*$, because it nontrivially involves the variable $k$. So the image of the complex elliptic genus

$$\varphi : MU_* \otimes Q \to \mathbf{Q}((y))[[q,k]]$$

is a graded polynomial ring $\mathbf{Q}[x_1, x_2, x_3, x_4]$, where $x_1$ is the image of $\mathbf{CP}^1$ and $x_2, x_3, x_4$ are the images of any three generators of the polynomial ring

$$MSU_* \otimes \mathbf{Q} = \mathbf{Q}[x_2, x_3, x_4, x_5, \ldots]$$

in complex dimensions 2, 3, 4, say those described earlier in this section.

## 4. Complex bordism modulo flops

Here is the first of the two main results of this paper.

THEOREM 4.1. *Let $I$ be the ideal in the complex bordism ring $MU_* \otimes \mathbf{Q}$ which is additively generated by differences $X_1 - X_2$, where $X_1$ and $X_2$ are smooth projective varieties related by a classical flop, as defined below. Then the complex elliptic genus, viewed as a ring homomorphism*

$$MU_* \otimes \mathbf{Q} \to \mathbf{Q}[x_1, x_2, x_3, x_4],$$

*is surjective with kernel equal to $I$. Equivalently, the Chern numbers (linear maps $MU_{2n} \otimes \mathbf{Q} \to \mathbf{Q}$) which are invariant under classical flops are exactly those which factor through the complex elliptic genus.*

After some geometric preliminaries, this result can be seen as a stronger form of one of Höhn's characterizations of the elliptic genus, in terms of twisted projective bundles [13, Satz 2.4.3].

Here two complex manifolds which are related by a classical flop, as defined below, are, in particular, both IH-small resolutions of the same complex space $Y$. So this theorem implies that the complex elliptic genus is an upper bound for the Chern numbers which can be defined for singular varieties. (Recall that in defining Chern numbers for singular varieties, we only allow definitions which are compatible with IH-small resolutions in the sense of Section 1.) We conjecture that this upper bound is an equality. Since the complex elliptic genus takes values in a proper quotient of

$$MU_* \otimes \mathbf{Q} = \mathbf{Q}[\mathbf{CP}^1, \mathbf{CP}^2, \ldots]$$

in complex dimensions $\geq 5$, we see in particular that not every Chern number can be defined for singular varieties in dimensions $\geq 5$.

The simplest singularity with two different IH-small resolutions is the 3-fold node $Y$ given in affine coordinates by $xy - zw = 0$, or equivalently, the cone over a smooth quadric surface $\mathbf{P}^1 \times \mathbf{P}^1 \subset \mathbf{P}^3$. Atiyah discussed the two small resolutions of this singularity [1]. Namely, blowing up the singular point of $Y$ gives a resolution $\tilde{X}$ which is not small; the inverse image of the singular point is a smooth divisor $\mathbf{P}^1 \times \mathbf{P}^1$ with normal bundle $O(-1,-1)$ (the tensor product of the line bundles $O(-1)$ on the two copies of $\mathbf{P}^1$). One can blow down either of the two families of $\mathbf{P}^1$'s on this divisor to give two resolutions of $Y$, $X_1$ and $X_2$, which are projective over $Y$, and both of which have fiber over the singular point of $Y$ equal to $\mathbf{P}^1$ with normal bundle $O(-1) \oplus O(-1)$. Now, if $Y$ is any projective 3-fold which is smooth outside one point, and $Y$ is Zariski locally isomorphic to the 3-fold node $xy - zw = 0$ near its singular point, then $Y$ clearly has two projective IH-small resolutions $X_1$ and $X_2$. We say that the smooth projective 3-folds $X_1$ and $X_2$ are related by a *classical flop.*

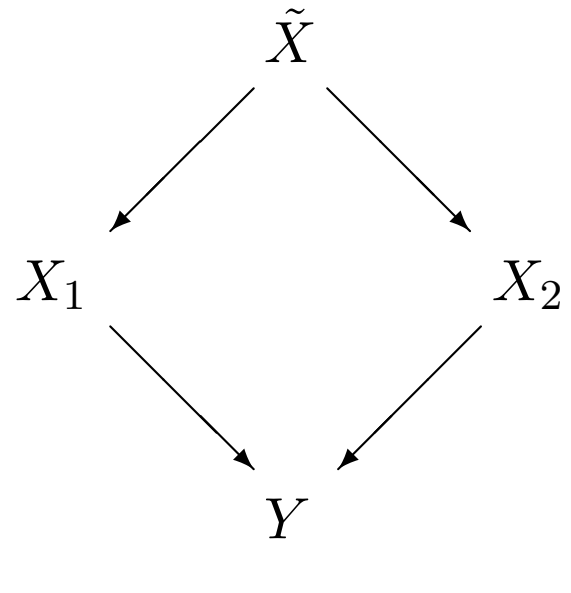

If we assume only that $Y$ is analytically locally isomorphic to $xy - zw = 0$, which is the usual definition of a node singularity of $Y$, then $Y$ still has two IH-small resolutions, but they need not be projective over $Y$, by [17, p. 171].

Our results, which are essentially topological, apply perfectly well to 3-fold nodes in this more general sense. The point of Theorem 4.1, however, is that even identifying manifolds related by a very special kind of flop reduces the complex bordism ring to the elliptic cohomology ring. It is to emphasize this point that we have defined "classical flops" in such a narrow sense.

To define classical flops in higher dimensions, let $Y$ be a singular projective $n$-fold which is Zariski locally isomorphic, near each point of its singular set $Z$, to the product of the 3-fold node with a smooth $(n-3)$-fold. Such a complex space $Y$ has two different IH-small resolutions $X_1$ and $X_2$ which are both smooth projective varieties. In this situation, we say that $X_1$ and $X_2$ are related by a classical flop. Geometrically, we are removing a $\mathbf{P}^1$-bundle over the smooth $(n-3)$-fold $Z$ inside $X_1$ and replacing it by a possibly different $\mathbf{P}^1$-bundle over $Z$.

*Proof of Theorem* 4.1. The first step is to show that the complex elliptic genus, applied to smooth projective varieties, is invariant under classical flops. This is a consequence of the crucial rigidity property of this genus, proved by both Krichever and Höhn. We need to state this not just for complex manifolds, but for weakly complex manifolds (those used in the definition of complex bordism). See Section 2 for the definitions of weakly complex manifolds and $SU$-manifolds.

Krichever-Höhn's rigidity theorem states that for any action of a compact connected Lie group $G$ on an $SU$-manifold $X$, the equivariant elliptic genus of $X$ is constant [19], [13]. The theorem has a consequence which we can state without mentioning equivariant genera. Namely, if $F \to E \to B$ is a fiber bundle of closed connected weakly complex manifolds, with structure group a compact connected Lie group $G$, and if $F$ is an $SU$-manifold, then the elliptic genus $\varphi$ satisfies

$$\varphi(E) = \varphi(F)\varphi(B),$$

by [13, Kor. 2.5.5]. (The condition on the structure group means that we start with a principal $G$-bundle over a weakly complex manifold $B$, and an action of $G$ on an $SU$-manifold $F$ which preserves the weakly complex structure, and then we let $E$ be the associated $F$-bundle over $B$.)

In fact, the elliptic genus $\varphi$, viewed as a surjection

$$\varphi : MU_* \otimes \mathbf{Q} \to \mathbf{Q}[x_1, x_2, x_3, x_4],$$

is the universal genus with the above multiplicativity property. Equivalently, the quotient of the complex cobordism ring $MU_* \otimes \mathbf{Q}$ by the relations $E - F \cdot B = 0$ for all fiber bundles as above is the polynomial ring $\mathbf{Q}[x_1, x_2, x_3, x_4]$. From this point of view, the remarkable fact about the elliptic genus is that this quotient ring is so big. If we divide out the complex cobordism ring by the relation $E = F \cdot B$ for all fiber bundles as above but without the $SU$ condition on $F$, then we get a much smaller and less interesting quotient ring of $MU_*\otimes\mathbf{Q}$, the image of the Hirzebruch $\chi_y$ genus. In other words, in some unexpected way, group actions on $SU$-manifolds are more restricted than group actions on general weakly complex manifolds. (The analogous statement which we get from the Landweber-Stong elliptic genus is that group actions on spin manifolds are more restricted than group actions on general oriented manifolds.)

We want to apply the multiplicativity property of the elliptic genus $\varphi$ to prove that $\varphi(X_1)=\varphi(X_2)$ for smooth projective $n$-folds $X_1$ and $X_2$ related by a classical flop, as defined above: $X_1$ and $X_2$ are the two IH-small resolutions of a singular projective variety $Y$ whose singular set $Z$ is a smooth subvariety of codimension 3, such that $Y$ is Zariski locally isomorphic near points of $Z$ to the product of an open subset of $Z$ with the 3-fold node. Thus $X_1$ and $X_2$ are isomorphic except over the inverse images of $Z$: the inverse image of $Z$ in $X_1$ is a $\mathbf{P}^1$-bundle $P(A)$ over $Z$, and the inverse image of $Z$ in $X_2$ is a $\mathbf{P}^1$-bundle $P(B)$ over $Z$ which may be different.

Since $X_1-P(A)$ is isomorphic to $X_2-P(B)$, the difference $X_1-X_2$ in the bordism group $MU_{2n}$ is equal to the class of a certain manifold $E$ which is fibered over $Z$. Namely, $E$ is the result of gluing a tubular neighborhood of $P(A)\subset X_1$ to a tubular neighborhood of $P(B)\subset X_2$ along their (diffeomorphic) boundaries. The manifold $E$ is not a complex manifold, but it is a weakly complex manifold in a natural way, with the given complex structure on the tubular neighborhood of $P(A)$ and with the negative weakly complex structure to the given one on the tubular neighborhood of $P(B)$. The negative of a weakly complex structure is discussed in Section 2, and the equality $X_1-X_2=E$ in bordism follows from Lemma 2.1.

The manifold $E$ is fibered over $Z$, with fiber a weakly complex 6-manifold $F$. One can construct $F$ as the difference of the two small resolutions of the 3-fold node, meaning the union of a tubular neighborhood of the $\mathbf{P}^1$ in one small resolution with a tubular neighborhood of the $\mathbf{P}^1$ in the other resolution along their common boundary, with the given complex structure on the first neighborhood and the negative weakly complex structure on the second one. This weakly complex 6-manifold $F$ was considered by Höhn [13, §1.3]; he called it the twisted projective space $\widetilde{\mathbf{CP}}_{2,2}$.

As a smooth manifold, $F$ is just $\mathbf{CP}^3$. But the crucial point is that the weakly complex structure defined here makes $F$ an $SU$-manifold. Indeed, $H^2(F,\mathbf{Z})=\mathbf{Z}$ maps isomorphically to $H^2$ of either of the two neighborhoods, so it is enough to check that $c_1(F)=0$ in one of those neighborhoods; but the first neighborhood is isomorphic to the bundle $O(-1)\oplus O(-1)$ over $\mathbf{P}^1$, which (as a 3-fold) has $c_1=0$.

The fiber bundle $F\to E\to Z$ has structure group $U(2)\times U(2)$. That is, there is an action of $U(2)\times U(2)$ on $F$, preserving the weakly complex structure, such that $E$ is the $F$-bundle over $Z$ associated to some $U(2)\times U(2)$- bundle over $Z$. Explicitly, given two rank-2 complex vector bundles $A$ and $B$ over $Z$, the weakly complex manifold $E$ is the "difference" between the bundle $B \otimes O(-1)$ over $P(A)$ and the bundle $A \otimes O(-1)$ over $P(B)$. Indeed, this is the way $E$ was constructed, since the inverse image of $Z$ in the smooth variety $X_1$ is isomorphic to a $\mathbf{P}^1$-bundle $P(A)$ over $Z$, with normal bundle of the form $B \otimes O(-1)$, and analogously for $X_2$. In the notation of Höhn [13, §1.3], $E$ is the twisted projective bundle $\widetilde{\mathbf{CP}}(A \oplus B)$.

Höhn describes the weakly complex structure on $E$, as follows. As a real manifold, $E$ is the $\mathbf{CP}^3$-bundle $P(A \oplus B^*)$ over $Z$. This manifold has a natural complex structure on the tangent bundle, and the weakly complex manifold $E$ will be defined by modifying this natural complex structure. We first need to describe the natural complex structure on the tangent bundle of $P(A \oplus B^*)$. Let $O(-1)$ denote the natural line subbundle of the pulled back rank-4 vector bundle $A \oplus B^*$ over $P(A \oplus B^*)$; then the tangent bundle along the fibers of $P(A \oplus B^*)$ is the tensor product of $(A \oplus B^*)/O(-1)$ with $O(1)$, that is, the quotient of $A \otimes O(1) \oplus B^* \otimes O(1)$ by a trivial bundle. So, as a $C^\infty$ complex vector bundle, the direct sum of the tangent bundle of $P(A \oplus B^*)$ with a trivial bundle is $A \otimes O(1) \oplus B^* \otimes O(1) \oplus TZ$.

Höhn defines the twisted projective bundle $E = \widetilde{\mathbf{CP}}(A \oplus B)$ to be the manifold $P(A \oplus B^*)$, with the complex structure on its stable tangent bundle changed to

$$T\widetilde{\mathbf{CP}}(A \oplus B) := A \otimes O(1) \oplus B \otimes O(-1) \oplus TZ.$$

(We are using here the fact that a complex vector bundle and its dual bundle can be identified as real vector bundles by choosing a hermitian metric on the given bundle.) It is elementary to identify this definition with the weakly complex structure on $E$ which comes from the above construction of $E$ by gluing.

Now we are in a position to apply Krichever-Höhn's rigidity theorem on the elliptic genus $\varphi$. Since $F \to E \to B$ is a fiber bundle of weakly complex manifolds with compact connected structure group $U(2) \times U(2)$ and $F$ is an $SU$-manifold, we have

$$\varphi(E) = \varphi(F)\varphi(B).$$

Moreover, the manifold $F$ is equal to 0 in the bordism group $MU_6$. The point is that for any 3-folds $X_1$ and $X_2$ related by a classical flop, we have $X_1 - X_2 = F$ in bordism. But flopping is a symmetric operation, so we also have $X_2 - X_1 = F$. Thus $2F = 0$ in bordism, and in fact $F = 0$ since the complex bordism ring is torsion-free. (We really only need that $F = 0$ in $MU_6 \otimes \mathbf{Q}$ for what follows.)

Since the elliptic genus is a homomorphism on the complex bordism ring, we have $\varphi(F) = 0$. Thus, by the above calculation, $\varphi(E) = 0$. Since we constructed $E$ as the difference in bordism between two complex $n$-folds $X_1$ and $X_2$ related by a classical flop, we have proved that $\varphi(X_1) = \varphi(X_2)$.

## 5. Complex bordism modulo flops, continued

We now finish the proof of Theorem 4.1. In the previous section, we showed that the elliptic genus is invariant under classical flops. It remains to show that the quotient of $MU_* \otimes \mathbf{Q}$ by differences $X_1 - X_2$, with $X_1$ and $X_2$ related by a classical flop, is no bigger than $\mathbf{Q}[x_1, x_2, x_3, x_4]$. There is a natural approach to proving this. Namely, in each complex dimension $n$ there is a Chern number $s_n$ such that an element $x$ of $MU_{2n} \otimes \mathbf{Q}$ is a polynomial generator of the ring $MU_* \otimes \mathbf{Q}$ if and only if $s_n(x)$ is not 0 [12]. Explicitly, $s_n$ is the $n^{\text{th}}$ power sum polynomial in the Chern classes $c_1, \ldots, c_n$; that is, as a polynomial in the Chern roots $x_1, \ldots, x_n$, $s_n$ is equal to $x_1^n + \cdots + x_n^n$. Our problem is solved if, for every $n \geq 5$, we can find complex $n$-folds $X_1$ and $X_2$ related by a classical flop such that $s_n(X_1) \neq s_n(X_2)$. Indeed, then the ideal of differences $X_1 - X_2$ contains a polynomial generator of $MU_* \otimes \mathbf{Q} = \mathbf{Q}[x_1, x_2, \ldots]$ in every degree at least 5, so the quotient by this ideal is at most $\mathbf{Q}[x_1, x_2, x_3, x_4]$, which is what we are trying to prove.

We repeat here that a classical flop is a diagram

$$\begin{array}{ccc} & \tilde{X} & \\ \swarrow & & \searrow \\ X_1 & & X_2 \\ \searrow & & \swarrow \\ & Y & \end{array} .$$

Here $Y$ is a singular projective $n$-fold which is Zariski locally isomorphic to the 3-fold node $xy - zw = 0$ times a smooth $(n-3)$-fold, near each point of its singular locus $Z$. We let $\tilde{X}$ be the blow-up of $Y$ along $Z$; $\tilde{X}$ is a smooth variety. The exceptional divisor $E \subset \tilde{X}$ is a $\mathbf{P}^1 \times \mathbf{P}^1$-bundle over the smooth $(n-3)$-fold $Z$. Finally, $X_1$ and $X_2$ are smooth varieties defined by contracting either of the two families of $\mathbf{P}^1$'s in $E \subset \tilde{X}$. In this situation, there are rank-2 vector bundles $A$ and $B$ on $Z$ such that the inverse image of $Z$ in $X_1$ is the $\mathbf{P}^1$-bundle $P(A)$ with normal bundle $B \otimes O(-1)$, and the inverse image of $Z$ in $X_2$ is $P(B)$ with normal bundle $A \otimes O(-1)$.

To compute the Chern number $s_n(X_1 - X_2)$, we could use Porteous's formula to relate the Chern numbers of $X_1$ and $X_2$ to those of their common blow-up $\tilde{X}$ [30]. It is more efficient, however, to use the observation in Section 4 that $X_1 - X_2$ is bordant to the twisted projective bundle $\widetilde{\mathbf{CP}}(A \oplus B)$ over $Z$. Notice that, for every smooth projective $(n-3)$-fold $Z$ with rank-2 algebraic vector bundles $A$ and $B$ over $Z$, there is a classical flop of $n$-folds in which $Z$, $A$, and $B$ play the roles explained above. That is, there is a singular projective $n$-fold $Y$ with singular set $Z$ such that $Y$ is, Zariski locally on $Z$, isomorphic to the 3-fold node $xy - zw = 0$ times an open subset of $Z$, and such that the blow-up $\tilde{X}$ of $Y$ along $Z$ has exceptional divisor $E = P(A) \times_Z P(B)$ with normal bundle $O(-1,-1) := O_{P(A)}(-1) \otimes O_{P(B)}(-1)$. Given $Z$, $A$, $B$, we can find such varieties $Y$, $\tilde{X}$, $X_1$, and $X_2$ as follows. First define a variety $E$ to be $P(A) \times_Z P(B)$, and then define $\tilde{X}$ to be the $\mathbf{P}^1$-bundle $P(O \oplus O(-1,-1))$ over $E$; then $E$ is embedded in $\tilde{X}$ in a natural way with normal bundle $O(-1,-1)$. We can then define $X_1$ and $X_2$ by contracting either of the two families of $\mathbf{P}^1$'s on $E \subset \tilde{X}$, and we can define $Y$ by contracting the family of $\mathbf{P}^1 \times \mathbf{P}^1$'s in $E \subset \tilde{X}$.

Thus Theorem 4.1 will be proved if we can find, for every $n \geq 5$, a smooth projective $(n-3)$-fold $Z$ with rank-2 algebraic vector bundles $A$ and $B$ such that $s_n(\widetilde{\mathbf{CP}}(A \oplus B)) \neq 0$.

We can compute $s_n(\widetilde{\mathbf{CP}}(A \oplus B))$ in slightly greater generality. For any weakly complex manifold $Z$ of real dimension $2(n-3)$ with $C^\infty$ rank-2 complex vector bundles $A$ and $B$, the definition of the twisted projective bundle $\widetilde{\mathbf{CP}}(A \oplus B)$ in Section 4 makes sense (as a weakly complex manifold). We will compute the Chern number $s_n$ of $E := \widetilde{\mathbf{CP}}(A \oplus B)$.

The Chern character of the tangent bundle of $E$ is

$$\begin{aligned} \mathrm{ch}(TE) &= \mathrm{ch}(A \otimes O(1) \oplus B \otimes O(-1) \oplus TZ) \\ &= e^u \mathrm{ch}(A) + e^{-u} \mathrm{ch}(B) + \mathrm{ch}(TZ), \end{aligned}$$

where $u$ denotes $c_1 O(1) \in H^2(E)$. If we write $x_1, x_2$ for the Chern roots of $A$ and $x_3, x_4$ for the Chern roots of $B$ (so, formally, $A$ has total Chern class $(1+x_1)(1+x_2)$ and $B$ has total Chern class $(1+x_3)(1+x_4)$), then the Chern character of $TE$ is

$$\begin{aligned} \mathrm{ch}(TE) &= e^u(e^{x_1} + e^{x_2}) + e^{-u}(e^{x_3} + e^{x_4}) + \mathrm{ch}(TZ) \\ &= e^{x_1+u} + e^{x_2+u} + e^{x_3-u} + e^{x_4-u} + \mathrm{ch}(TZ). \end{aligned}$$

Since $s_n(TE) = n!\,\mathrm{ch}_n(TE)$, it follows that

$$s_n(TE) = (x_1+u)^n + (x_2+u)^n + (x_3-u)^n + (x_4-u)^n + s_n(TZ).$$

Here $s_n(TZ) = 0$ since the bundle $TZ$ on $E$ is pulled back from $Z$, which has real dimension $2(n-3)$. So in fact

$$s_n(TE) = (x_1+u)^n + (x_2+u)^n + (x_3-u)^n + (x_4-u)^n.$$

To compute $s_nE := \int_E s_n(TE)$, we rewrite this integral as an integral over $Z$, by the equality

$$s_nE = \int_Z \pi_* s_n(TE),$$

where $\pi : E \to Z$ is the projection and $\pi_* : H^iE \to H^{i-6}Z$ is the corresponding pushforward map. Here $E$ is identified with the $\mathbf{CP}^3$-bundle $P(A \oplus B^*)$ as a real manifold, and one checks that this identification preserves orientation, so it is enough to describe the pushforward map for complex projective bundles. Namely, $H^*E$ is a free module over $H^*Z$ with basis $1, u, u^2, u^3$, and the pushforward map $\pi_*$ is a map of $H^*Z$-modules, so it is enough to describe $\pi_*$ on the basis elements, which is easy: $\pi_*(u^i) = 0$ for $0 \leq i \leq 2$ and $\pi_*(u^3) = \int_{\mathbf{CP}^3} u^3 = 1$. In fact, to apply $\pi_*$ to $s_nTE$ as computed above, it is convenient to have a formula for $\pi_*(u^i)$ for any $i \geq 0$. These are given by Segré classes, that is, inverse Chern classes.

LEMMA 5.1. *For any space $X$ with a complex vector bundle $V$ of rank $r$ over $X$, let $\pi : P(V) \to X$ be the projective bundle of lines in $V$ and let $u = c_1O(1) \in H^2P(V)$. Then*

$$\pi_*(u^i) = c_{i-(r-1)}(-V)$$

*for all $i \geq 0$.*

A reference for the lemma is Fulton's book [9, p. 47]. Thus, for the bundle $\pi : E \to Z$, we have

$$\pi_*(u^i) = c_{i-3}(-(A \oplus B^*)).$$

In terms of the Chern roots $x_1, x_2$ of $A$ and $x_3, x_4$ of $B$, we have

$$\begin{aligned}\pi_*(u^i) &= \sum_{\substack{i_1+i_2+i_3+i_4=i-3\\ i_j\geq 0}} (-x_1)^{i_1}(-x_2)^{i_2}x_3^{i_3}x_4^{i_4}\\ &= \sum_{\substack{i_1+i_2+i_3+i_4=i-3\\ i_j\geq 0}} (-1)^{i_1+i_2}x_1^{i_1}x_2^{i_2}x_3^{i_3}x_4^{i_4}.\end{aligned}$$

So we get the following expression for the Chern number $s_n$ of $E$:

$$\begin{aligned}s_nE &= \int_Z \pi_* s_nTE\\ &= \int_Z \pi_*\Big[(x_1+u)^n + (x_2+u)^n + (x_3-u)^n + (x_4-u)^n\Big]\\ &= \int_Z \sum_{i=3}^n \binom{n}{i}\Big[x_1^{n-i} + x_2^{n-i} + (-1)^i x_3^{n-i} + (-1)^i x_4^{n-i}\Big]\pi_* u^i.\end{aligned}$$

(The sum is written only over $i \geq 3$ because $\pi_*(u^i) = 0$ for $i \leq 2$.) When we plug in the formula for $\pi_*(u^i)$ and use the identity $\sum_{j=0}^{i}(-1)^j\binom{n}{j} = (-1)^i\binom{n-1}{i}$, we get the definitive formula for $s_nE$:

$$\begin{aligned} s_nE &= \int_Z \sum_{\substack{i_1+i_2+i_3+i_4=n-3\\ i_r\geq 0}} x_1^{i_1}x_2^{i_2}x_3^{i_3}x_4^{i_4}\Bigg[(-1)^{i_2}\binom{n-1}{i_1}+(-1)^{i_1}\binom{n-1}{i_2} \\ &\qquad\qquad + (-1)^{i_4+1}\binom{n-1}{i_3}+(-1)^{i_3+1}\binom{n-1}{i_4}\Bigg]. \end{aligned}$$

Using this formula, it is easy to find, for every $n \geq 5$, a smooth projective $(n-3)$-fold $Z$ with rank-2 algebraic vector bundles $A$ and $B$ such that the associated twisted projective bundle $E$ has $s_n(E) \neq 0$. As explained earlier in this section, this will complete the proof of Theorem 4.1. Let $Z = \mathbf{CP}^{n-3}$, with $A = O(1) + O$ and $B = O^{\oplus 2}$. Then, in the above notation, $x_1 = c_1O(1)$ on $\mathbf{CP}^{n-3}$ and $x_2 = x_3 = x_4 = 0$, and so

$$\begin{aligned} s_nE &= \int_Z x_1^{n-3}\left[\binom{n-1}{n-3}+(-1)^{n-3}\binom{n-1}{0}-\binom{n-1}{0}-\binom{n-1}{0}\right] \\ &= (n^2-3n+2(-1)^{n+1}-2)/2 \\ &= \begin{cases} n(n-3)/2 & \text{if } n \text{ is odd} \\ (n+1)(n-4)/2 & \text{if } n \text{ is even.} \end{cases} \end{aligned}$$

Thus $s_nE$ is 0 for $n = 3$ and $n = 4$, but nonzero for all $n \geq 5$ (Theorem 4.1). ∎

## 6. *SU*-bordism modulo flops

In this section, we give a geometric description of the kernel of the complex elliptic genus restricted to $MSU_* \otimes \mathbf{Z}[1/2]$. Namely, this kernel is equal to the ideal $I$ in $MSU_* \otimes \mathbf{Z}[1/2]$ generated by twisted projective bundles $\widetilde{\mathbf{CP}}(A \oplus B)$ over weakly complex manifolds $Z$ such that the complex vector bundles $A$ and $B$ over $Z$ have rank 2 and $c_1Z + c_1A + c_1B = 0$; in this case, the total space is an $SU$-manifold. In view of Section 4, it is reasonable to call $I$ the ideal of $SU$-flops: some elements of $I$ will arise geometrically from birational equivalences between compact complex manifolds with trivial canonical bundle. At the same time, we find that $MSU_* \otimes \mathbf{Z}[1/2]/I$ is a polynomial ring of the form $\mathbf{Z}[1/2][x_2, x_3, x_4]$. These results are analogous to the results of Kreck and Stolz, describing the kernel of the Ochanine genus on $M\mathrm{Spin}_*$ in terms of $\mathbf{HP}^2$-bundles, except that for now we work away from the prime 2 [18].

*Remarks.* (1) For integral questions such as this, it seems more natural to work with the ring $MSU_*$ rather than $MU_*$, for example because the image of $MU_*$ under the complex elliptic genus is not finitely generated, although after tensoring with $\mathbf{Q}$ it becomes the polynomial ring $\mathbf{Q}[x_1, x_2, x_3, x_4]$. In fact, even the image of the $\chi_y$ genus on $MU_*$ is not finitely generated.

The image of $MSU_*$, on the other hand, is quite simple, as explained above. In particular, given the results above, the Sullivan-Baas method of bordism with singularities produces a multiplicative cohomology theory which is a module over $MSU \otimes \mathbf{Z}[1/2]$ and which has coefficient ring $\mathbf{Z}[1/2][x_2, x_3, x_4]$ [2], [24]. (The Sullivan-Baas method gives a cohomology theory with the coefficient ring we want because the ideal $I$ is defined by a regular sequence in the ring

$$MSU_* \otimes \mathbf{Z}[1/2] = \mathbf{Z}[1/2][x_2, x_3, x_4, \ldots],$$

as follows from the above results.) This is a natural integral version of complex elliptic cohomology theory (at least over $\mathbf{Z}[1/2]$), defined here for the first time.

(2) One might define an ideal $I \subset MSU_*$ of "$SU$-flops" in several other ways. In particular, it would be closer to Kreck and Stolz's description of the kernel of the Ochanine genus via $\mathbf{HP^2}$-bundles to consider only twisted projective bundles whose base as well as whose total space is an $SU$-manifold [18]. This seems to be definitely the wrong definition at the prime 2, so we have preferred the more general definition of $SU$-flops above, which at least has a chance of giving the "right" ideal in $MSU_*$ as well as in $MSU_* \otimes \mathbf{Z}[1/2]$. At the prime 2, it may also be necessary to consider twisted projective bundles with structure group $(U(2) \times U(2))/U(1)$ rather than $U(2) \times U(2)$.

We now turn to the proof of this section's theorem:

THEOREM 6.1. *The kernel of the complex elliptic genus on* $MSU_* \otimes \mathbf{Z}[1/2]$ *is equal to the ideal* $I$ *of SU-flops, as defined above. Also, the quotient ring is a polynomial ring*:

$$MSU_* \otimes \mathbf{Z}[1/2]/I \cong \mathbf{Z}[1/2][x_2, x_3, x_4].$$

*Proof.* We use Novikov's description of the ring $MSU_* \otimes \mathbf{Z}[1/2]$ [29]. It is a graded polynomial ring

$$MSU_* \otimes \mathbf{Z}[1/2] = \mathbf{Z}[1/2][x_2, x_3, x_4, \ldots],$$

$x_n \in MSU_{2n}$. An $SU$-manifold $X$ of real dimension $2n$, $n \geq 2$, is a polynomial generator of $MSU_* \otimes \mathbf{Z}[1/2]$ if and only if

$$s_n X = \begin{cases} \pm p(\text{a power of } 2) & \text{if } n \text{ is a power of an odd prime } p, \\ \pm p(\text{a power of } 2) & \text{if } n+1 \text{ is a power of an odd prime } p, \\ \pm(\text{a power of } 2) & \text{otherwise.} \end{cases}$$

We will show that for every $n \geq 5$, the greatest common divisor of the integers $s_n X$ for $SU$-flops $X$ of real dimension $2n$ is as above. This will imply the theorem, as follows. The statement means that the ideal $I \subset MSU_* \otimes$ $\mathbf{Z}[1/2]$ contains a polynomial generator of $MSU_* \otimes \mathbf{Z}[1/2]$ in real dimension $2n$ for all $n \geq 5$. So $MSU_* \otimes \mathbf{Z}[1/2]/I$ is a quotient of the ring $\mathbf{Z}[1/2][x_2, x_3, x_4]$. But the complex elliptic genus gives a homomorphism of graded rings

$$MSU_* \otimes \mathbf{Z}[1/2]/I \to \mathbf{Q}[x_2, x_3, x_4],$$

where we know that $I$ maps to 0 since this genus is 0 on all flops, by Theorem 4.1. Moreover, this homomorphism is surjective after tensoring with $\mathbf{Q}$ [13]. Since $MSU_* \otimes \mathbf{Q}/I$ is a quotient of $\mathbf{Q}[x_2, x_3, x_4]$, a consideration of dimensions shows that $MSU_* \otimes \mathbf{Q}/I$ is actually equal to $\mathbf{Q}[x_2, x_3, x_4]$. Since $MSU_* \otimes \mathbf{Z}[1/2]/I$ is a quotient of the torsion-free ring $\mathbf{Z}[1/2][x_2, x_3, x_4]$, any relation would show up rationally, and so $MSU_* \otimes \mathbf{Z}[1/2]/I$ is equal to $\mathbf{Z}[1/2][x_2, x_3, x_4]$. Finally, because this ring is torsion-free, the complex elliptic genus

$$MSU_* \otimes \mathbf{Z}[1/2]/I \to \mathbf{Q}[x_2, x_3, x_4]$$

is injective, since this is true rationally. That is, $I$ is exactly the kernel of the complex elliptic genus on $MSU_* \otimes \mathbf{Z}[1/2]$.

Thus the theorem will be proved if we can show that the greatest common divisor of the integers $s_n X$ for $SU$-flops $X$ is as promised above. Our tool will be the following lemma.

LEMMA 6.2. *For a weakly complex manifold $Z$ of real dimension $2n$ with complex line bundles $L_1$, $L_2$, $L_3$ on $Z$, consider the set of integers*

$$c_1^{i_1} L_1 c_1^{i_2} L_2 c_1^{i_3} L_3 c_1^{i_4} Z,$$

*for all natural numbers $i_1, \ldots, i_4$ with $i_1 + i_2 + i_3 + i_4 = n$. This gives a homomorphism $f : MU_{2n}(BU(1)^3) \to \mathbf{Z}^N$, where $N$ is the number of partitions $i_1 + \cdots + i_4 = n$. Then $f$ becomes surjective after tensoring with $\mathbf{Z}[1/2]$.*

*Proof of Lemma* 6.2. The point is that for every $n \geq 0$, there is a weakly complex manifold $F$ of real dimension $2n$ with $c_1^n(F)$ a unit in $\mathbf{Z}[1/2]$. (We cannot always make $c_1^n(F)$ a unit in $\mathbf{Z}$; for example, $c_1$ of every complex curve is even.) For $n \geq 1$, we can take $F$ to be a suitable $\mathbf{Z}$-linear combination of $\mathbf{P}^n$ and $\mathbf{P}^1 \times \mathbf{P}^{n-1}$, since

$$c_1^n(\mathbf{P}^n) = (n+1)^n$$

and

$$c_1^n(\mathbf{P}^1 \times \mathbf{P}^{n-1}) = 2n^n.$$

Also, we use the fact that for every $n \geq 0$, the map $MU_{2n}(BU(1)^3) \to H_{2n}(BU(1)^3, \mathbf{Z})$ is surjective, as is true for any space with torsion-free cohomology in place of $BU(1)^3$. (This follows from inspection of the Atiyah-Hirzebruch spectral sequence for bordism.) Equivalently, for each $n \geq 0$, there are weakly complex manifolds $M$ of real dimension $2n$ with complex line bundles $L_1, L_2, L_3$ such that the integers $c_1^{i_1}L_1c_1^{i_2}L_2c_1^{i_3}L_3$, $i_1+i_2+i_3=n$, are whatever we like.

The lemma follows from combination of these two observations. For $0 \leq k \leq n$, consider weakly complex manifolds of the form $X = F \times M$, where $F$ is a fixed manifold of real dimension $2k$ such that $c_1^k(F)$ is a unit in $\mathbf{Z}[1/2]$, and $M$ has real dimension $2(n-k)$ and has three complex line bundles $L_1, L_2, L_3$ on it. We can view $L_1, L_2, L_3$ as line bundles on the product manifold $X$. The Chern numbers in $c_1X, c_1L_1, c_1L_2, c_1L_3$ of degree $< k$ in $c_1X$ are 0, and those of degree $k$ in $c_1X$ are equal to $c_1^k(F)$ times the Chern numbers in $c_1L_1, c_1L_2, c_1L_3$ on $M$. By taking a suitable $\mathbf{Z}[1/2]$-linear combination of these manifolds $X$ over $0 \leq k \leq n$, we can make the Chern numbers in $c_1X, c_1L_1, c_1L_2, c_1L_3$ anything we like over $\mathbf{Z}[1/2]$. □

Using Lemma 6.2, we proceed to prove Theorem 6.1. For any weakly complex manifold $Z$ of real dimension $2(n-3)$ with complex line bundles $L_1, L_2, L_3$ on $Z$, we define rank-2 bundles $A$ and $B$ on $Z$ by $A = L_1 \oplus L_2$ and $B = L_3 \oplus (K_Z \otimes L_1^* \otimes L_2^* \otimes L_3^*)$. These are chosen so that $c_1Z + c_1A + c_1B = 0$, which is the condition needed to ensure that the twisted projective bundle $E := \widetilde{\mathbf{CP}}(A \oplus B)$ is an $SU$-manifold, that is, by our definition, an $SU$-flop.

By Section 5, the Chern number $s_nE$ is a certain explicit linear combination of the integers

$$c_1^{i_1}L_1c_1^{i_2}L_2c_1^{i_3}L_3c_1^{i_4}Z.$$

By Lemma 6.2, the greatest common divisor of the integers $s_nE$ obtained this way, in the ring $\mathbf{Z}[1/2]$, is simply the greatest common divisor of the coefficients of the integers $c_1^{i_1}L_1c_1^{i_2}L_2c_1^{i_3}L_3c_1^{i_4}Z$ in the formula for $s_nE$. By inspection of that formula, this greatest common divisor is equal to the greatest common divisor of the integers

$$(-1)^{i_2}\binom{n-1}{i_1} + (-1)^{i_1}\binom{n-1}{i_2} + (-1)^{i_4+1}\binom{n-1}{i_3} + (-1)^{i_3+1}\binom{n-1}{i_4}$$

over all partitions $i_1 + \cdots + i_4 = n-3$. Thus Theorem 6.1 will follow if we can show, for $n \geq 5$, that if an odd prime number $p$ divides all these integers, then either $n$ or $n+1$ is a power of $p$, and in those cases one of these integers is not divisible by $p^2$.

We will only use a few of these integers, the ones with $i_2 = i_4 = 0$. Write $i$ for $i_1$, so that $0 \leq i \leq n-3$ and $i_3 = n-3-i$. Then the above integer is

$$\begin{aligned} &= \binom{n-1}{i} + (-1)^i - \binom{n-1}{n-3-i} + (-1)^{n-3-i+1} \\ &= \binom{n}{i+1} - \binom{n}{i+2} + (-1)^i[1+(-1)^n], \end{aligned}$$

by the identity $\binom{n}{i} = \binom{n-1}{i} + \binom{n-1}{i-1}$. Thus, if $n$ is odd, the above integer is $\binom{n}{i+1} - \binom{n}{i+2}$, for $0 \le i \le n-3$. If $n$ is even, since we are interested in the greatest common divisor of these numbers, we can take the sum of the above number for $i$ and the above number for $i-1$; this gives the integers $\binom{n+1}{i+1} - \binom{n+1}{i+2}$, for $1 \le i \le n-3$, and (still for $n$ even) we can also remember the above number for $i = 0$, which is $2\binom{n+1}{1} - \binom{n+1}{2}$.

Suppose first that $n$ is odd (and, as always, $n \ge 5$). We will show that if an odd prime number $p$ divides $\binom{n}{i+1} - \binom{n}{i+2}$ for all $0 \le i \le n-3$, then $n$ is a power of $p$, and one of these numbers is not divisible by $p^2$; this will prove what we want for $n$ odd. Since $n \ge 5$, our assumption implies in particular that $p$ divides $\binom{n}{2} - \binom{n}{1} = n(n-3)/2$ and $\binom{n}{3} - \binom{n}{2} = n(n-1)(n-5)/6$. If $p$ does not divide $n$, then $n \equiv 3 \pmod p$ and $n \equiv 1$ or $5 \pmod p$, a contradiction. So $p$ does divide $n$.

Equivalently, $p$ divides $\binom{n}{1}$; so our assumption tells us that $p$ divides $\binom{n}{i}$ for all $1 \le i \le n-1$. It is elementary that this implies that $n$ is a power of $p$. Also, it is then elementary that the greatest common divisor of the integers $\binom{n}{i}$ is exactly $p$.

It remains to show that one of the differences $\binom{n}{i+1} - \binom{n}{i+2}$ is not divisible by $p^2$, $0 \le i \le n-3$. If they are all zero modulo $p^2$, then since one of the integers $\binom{n}{i}$, $1 \le i \le n-1$, is nonzero modulo $p^2$, they are all nonzero modulo $p^2$. In particular this applies to $\binom{n}{1} = n$, which is a power of $p$; so $n = p$. In this case, we have $p(p-3)/2 = \binom{n}{2} - \binom{n}{1} \equiv 0 \pmod{p^2}$, and so $n = p = 3$. This contradicts our assumption that $n \ge 5$. So one of the differences $\binom{n}{i+1} - \binom{n}{i+2}$ is not divisible by $p^2$, and our proof is complete for $n$ odd.

It remains to consider $n$ even (with, as always, $n \ge 5$). We will show that if an odd prime number $p$ divides $\binom{n+1}{i+1} - \binom{n+1}{i+2}$ for all $1 \le i \le n-3$ as well as $2\binom{n+1}{1} - \binom{n+1}{2}$, then $n+1$ is a power of $p$, and one of the integers mentioned is not divisible by $p^2$. This will complete the proof of Theorem 6.1.

Since $p$ divides $2\binom{n+1}{1} - \binom{n+2}{2} = -(n+1)(n-4)/2$, $n+1$ is congruent to 0 or 5 modulo $p$. Suppose that $n+1$ is not a multiple of $p$, so that $n+1 \equiv 5 \pmod p$. Since $n \ge 5$, $n+1$ is at least $p+5$ in this case. Our assumptions imply, in this case, that the numbers $\binom{n+1}{i}$ for $2 \le i \le n-1$ are all equal and nonzero modulo $p$. But

$$\binom{n+1}{p+1} = \binom{n+1}{p} \frac{n+1-p}{p+1},$$

so that the fraction on the right must equal $1 \pmod p$, which says that $n+1 \equiv 1 \pmod p$, contradicting our assumption that $n+1 \equiv 5 \pmod p$. Now, in fact, $n+1$ is a multiple of $p$. By our assumptions, then, the numbers $\binom{n+1}{i}$ for $1 \le i \le n$ are all $0 \pmod p$. It follows that $n+1$ is a power of $p$.

Then the greatest common divisor of the integers $\binom{n+1}{i}$, $1 \le i \le n$, is $p$. We have to show that one of the differences $\binom{n+1}{i+1} - \binom{n+1}{i+2}$, $1 \le i \le n-3$, or else the number $2\binom{n+1}{1}-\binom{n+1}{2}$, is not divisible by $p^2$. If these numbers are all 0 modulo $p^2$, then all the numbers $\binom{n+1}{i}$, $1\leq i\leq n$, are nonzero modulo $p^2$ since one of them is. In particular this applies to $n+1$, which is a power of $p$, and so $n+1=p$. But then $p^2$ divides $2\binom{n+1}{1}-\binom{n+1}{2}=-p(p-5)/2$, so that $n+1=p=5$, contradicting our assumption that $n\geq 5$. So we have proved that one of the numbers $\binom{n+1}{i+1}-\binom{n+1}{i+2}$, $1\leq i\leq n-3$, or else the number $2\binom{n+1}{1}-\binom{n+1}{2}$, is not divisible by $p^2$. This completes the proof of Theorem 6.1. □

## 7. Saito's homology classes $\chi_i^{n-k}$

In this section we explain how Morihiko Saito's definition of a pure Hodge structure on intersection homology [32], [33] implicitly includes a definition of certain natural homology classes $\chi_i^{n-k}$ on a singular algebraic variety. We will use these classes to define some new Chern numbers for singular varieties in Section 8.

Saito showed that for any complex algebraic variety $X$, the intersection homology complex $IC_X$ in the derived category of $\mathbf{C}_X$-modules has a natural "Hodge" filtration $F$ in the derived category such that the associated graded objects $\mathrm{Gr}_p^F IC_X$ live naturally in the derived category $D^b_{\mathrm{coh}}(O_X)$ of bounded complexes of $O_X$-modules with cohomology sheaves which are coherent $O_X$-modules [33, p. 273]. If $X$ is smooth, so that $IC_X=\mathbf{C}_X\in D(\mathbf{C}_X)$, $IC_X$ is quasi-isomorphic to the de Rham complex

$$0\to\Omega_X^0\to\Omega_X^1\to\cdots,$$

and the filtration $F$ is the obvious filtration so that $\mathrm{Gr}_p^F IC_X=\Omega_X^p[p]$. Thus, for a general variety $X$, the object $\mathrm{Gr}_p^F IC_X\in D^b_{\mathrm{coh}}(O_X)$ is a generalization of the sheaf of $p$-forms on a smooth variety. (A different generalization of the sheaf of $p$-forms on a smooth variety to an object in the derived category was found earlier by du Bois [7], related to ordinary cohomology rather than intersection cohomology. Saito was partly inspired by du Bois's work.)

The filtration $F$ is preserved under proper pushforward in a precise sense [33, p. 273]. In particular, if $f:X\to Y$ is an IH-small resolution, then $f_*IC_X=IC_Y\in D(C_X)$ compatibly with the filtrations $F$ on $IC_X$ and $IC_Y$, and $f_*\mathrm{Gr}_p^F IC_X=\mathrm{Gr}_p^F IC_Y\in D^b_{\mathrm{coh}}(O_Y)$. Here $f_*$ means $Rf_*$, as is natural in derived categories.

In particular, the alternating sum

$$\chi_p:=\sum_i(-1)^{i+p}\mathcal{H}^i\mathrm{Gr}_p^F IC_X$$

is an element of the Grothendieck group $G_0X$ of coherent sheaves on any variety $X$. It is equal to the class of $\Omega^p_X$ for $X$ smooth, and it satisfies $f_*\chi_p(X) = \chi_p(Y)$ for an IH-small resolution $f : X \to Y$.

Finally, we can apply Baum-Fulton-MacPherson's natural homomorphism [4] from the Grothendieck group $G_0X$ of coherent sheaves to topological $K$-homology $K_0^{\text{top}}X$, followed by the homological version of the Chern character $\text{hh} : K_0^{\text{top}}X \to H_*(X, \mathbf{Q})$ (also described in [4]), to $\chi_p(X) \in G_0X$. Let $\chi^r_p(X)$ be the part of $\text{hh}(\chi_p(X))$ in $H_{2(n-r)}(X, \mathbf{Q})$, where $n = \dim X$. Then $\chi^{n-k}_p$ is a homology class naturally associated to any complex algebraic variety $X$ (it lives in Borel-Moore homology if $X$ is noncompact). If $X$ is smooth, it is the degree $n-k$ part of the cohomology class $\text{td}(T_X)\text{ch}(\Omega^p_X)$, since the homology Chern character hh takes a vector bundle $E$ to the Chern character of $E$ times the Todd class of $X$. Finally, the homology Chern character is a natural transformation on the category of proper algebraic maps; so if $f : X \to Y$ is an IH-small resolution, then $f_*\chi^{n-k}_p(X) = \chi^{n-k}_p(Y)$.

## 8. The Chern numbers $c_1^k\chi_i^{n-k}$

Our positive results about the Chern numbers $c_1^k\chi_i^{n-k}$ for $n$-folds come in two slightly different forms, both to be proved in this section. (We describe these Chern numbers in terms of a certain genus, the $\chi_{yz}$ genus, in Section 9.) First, we have:

THEOREM 8.1. *For any singular projective variety $Y$ with two projective IH-small resolutions $X_1$ and $X_2$,*

$$c_1^k\chi_i^{n-k}(X_1) = c_1^k\chi_i^{n-k}(X_2)$$

*for all $0 \le k \le n$, $0 \le i \le n-k$.*

The projectivity assumptions can be weakened, as follows: it suffices to let $Y$ be a compact complex space with two IH-small resolutions $X_i \to Y$, $i = 1, 2$, which are projective over an analytic neighborhood of each point of $Y$. In fact, even this weaker projectivity assumption should be irrelevant.

An equivalent statement is that there exists some extension of the Chern number $c_1^k\chi_i^{n-k}$ to singular $n$-folds which agrees with the corresponding number for any projective IH-small resolution. Of course, it would be more satisfying to define explicitly at least one extension of the Chern number $c_1^k\chi_i^{n-k}$ to singular $n$-folds, and we can do this for those varieties which have a relative canonical model; again, the minimal model conjecture (Conjecture 0-4-4 in [14]) would imply that every variety has a relative canonical model (Theorem 3-3-1 in [14]).

We recall some of the relevant definitions. A relative canonical model for a variety $Y$ is defined, starting from any resolution of singularities $f : X \to Y$, as the variety $X_0 := \mathrm{Proj}(\oplus_{n\geq 0} f_* K_X^{\otimes n}) \to Y$, assuming that the sheaf of $O_Y$-algebras $\oplus_{n\geq 0} f_* K_X^{\otimes n}$ on $Y$ is locally finitely generated [14, p. 301]. The sheaf of algebras $\oplus_{n\geq 0} f_* K_X^{\otimes n}$ on $Y$ is independent of the resolution $X$ of $Y$; this is the classical fact that sections of pluricanonical bundles on smooth varieties are birationally invariant. So every variety has at most one relative canonical model. This makes the relative canonical model a useful tool for defining invariants of singular varieties whenever it can be shown to exist. It exists for varieties of dimension at most 3 [26, Th. 0.3.12], and for varieties with toroidal singularities. (Reid [31, Th. 0.2] proves that it exists for global toric varieties, and the definition of the relative canonical model is analytically local on $Y$.)

A crucial property of relative canonical models is that, as Reid found, their singularities are not too bad. In particular, if $X_0$ is the relative canonical model of any variety $Y$, then $X_0$ has canonical singularities (see [17, p. 121] for the definition), and $X_0$ is $\mathbf{Q}$-Gorenstein, which means that some power of the canonical class $K_{X_0}$ is a line bundle on $X_0$. Thus, we have a well-defined cohomology class $c_1(X_0) := -K_{X_0} \in H^2(X_0, \mathbf{Q})$.

Now we can state the second version of our positive result on the Chern numbers $c_1^k \chi_i^{n-k}$.

*Definition.* Let $Y$ be a singular projective variety which has a relative canonical model $X_0$. Then we define

$$c_1^k \chi_i^{n-k}(Y) := c_1^k \chi_i^{n-k}(X_0) \in \mathbf{Q}.$$

The formula on the canonical model $X_0$ makes sense because $\chi_i^{n-k}$ can be defined as a homology class on arbitrary varieties (see Section 7), and $c_1 = -K_{X_0}$ is a cohomology class in $H^2(X_0, \mathbf{Q})$ since the canonical model $X_0$ is $\mathbf{Q}$-Gorenstein. The formula gives a well-defined rational number associated to $Y$ because the canonical model is unique. However, we also want to know that this definition is compatible with IH-small resolutions in the sense required, and that turns out to be true:

THEOREM 8.2. *Let $Y$ be a singular projective variety which has a projective* IH-*small resolution $X \to Y$. Then $Y$ also has a relative canonical model $X_0$, which factors $X \to X_0 \to Y$, and*

$$c_1^k \chi_i^{n-k}(X) = c_1^k \chi_i^{n-k}(X_0).$$

This theorem will be deduced from the fact, of interest in its own right, that projective IH-small resolutions are relative minimal models in the sense of Mori's theory, as we will explain. This fact is a restatement of a theorem of Wisniewski's [34].

Now we turn to the proofs.

*Proof of Theorem* 8.1. This will follow from Theorem 8.2. Indeed, if $X_1$ and $X_2$ are projective IH-small resolutions of $Y$, then by Theorem 8.2, $Y$ has a relative canonical model $X_0$, lying under both $X_1$ and $X_2$,

$$\begin{array}{ccc} X_1 & & X_2 \\ & \searrow \quad \swarrow & \\ & X_0 & \\ & \downarrow & \\ & Y, & \end{array}$$

and $c_1^k\chi_i^{n-k}(X_1) = c_1^k\chi_i^{n-k}(X_0) = c_1^k\chi_i^{n-k}(X_2)$. □

*Proof of Theorem* 8.2. As mentioned above, the main point is the following proposition, a restatement of a theorem of Wisniewski's.

PROPOSITION 8.3. *An* IH-*small resolution* $f : X \to Y$ *such that* $X$ *is projective over* $Y$ *is a relative minimal model.*

*Proof of Proposition* 8.3. By definition, a relative minimal model of a variety $Y$ is a variety $X$ with $\mathbf{Q}$-factorial terminal singularities together with a projective birational morphism $f : X \to Y$ such that the canonical class $K_X$ is $f$-nef; that is, $K_X \cdot C \geq 0$ for all curves $C$ in $X$ which map to a point in $Y$. Again, some references are [31, Th. 0.2] for the toric case, [17, p. 121] for the definition of terminal singularities, and [14] for the most detailed development of the theory. Kollár [16] gives a good introductory survey, and his more recent survey [17] is also very useful.

If $X$ is a projective IH-small resolution of a variety $Y$, then $X$ is smooth and hence has $\mathbf{Q}$-factorial terminal singularities. So to show that $X$ is a relative minimal model we only have to show that $K_X$ is $f$-nef. Suppose that $K_X$ is not $f$-nef. Let $N(X/Y)$ be the real vector space spanned by the curves on $X$ which map to a point in $Y$ modulo numerical equivalence, that is, modulo the relation that a linear combination of curves is 0 if it has 0 intersection number with every line bundle on $X$. Define $\overline{NE(X/Y)}$, the cone of curves in $X$ over $Y$, to be the closed cone in $N(X/Y)$ generated by curves in $X$ which map to a point in $Y$. In these terms, since $K_X$ is not $f$-nef, the cone $\overline{NE(X/Y)} \cap \{z \in N(X/Y) : K_X \cdot z < 0\}$ is nonempty. By Mori [25, Th. 1.4], later generalized by several people (see [14, Ch. 4]), this intersection is locally a rational polyhedral cone, so in particular, given that it is nonempty, it has an extremal ray $R$; let us pick one. By Mori, Kawamata, Benveniste, Reid, Ando, and Shokurov (see [14, Ch. 3]), we can contract the extremal ray $R$. This means that there is a normal variety $X'$ with surjective morphisms $X \to X' \to Y$ such that $X'$ is projective over $Y$, and such that a curve $C$ in $X$ which maps to a point in $Y$ also maps to a point in $X'$ if and only if the class of $C$ in $N(X/Y)$ lies on the ray $R$. The variety $X'$ is uniquely defined by these properties. Since the map $X \to Y$ is IH-small, so is the contraction map $X \to X'$. We will derive a contradiction from the existence of such an IH-small contraction, thus proving that $X$ was in fact a relative minimal model of $Y$.

Define the *length* of the extremal ray $R$ to be

$$l(R) := \min \{-K_X \cdot C : [C] \in R - \{0\}, C \text{ a rational curve}\};$$

thus $l(R)$ is a positive integer, given that $X$ is smooth so that $-K_X \cdot C$ is always an integer. Now we can state Wisniewski's result [34, Th. 1.1].

THEOREM 8.4. *Let $R$ be an extremal ray in a smooth variety $X$, giving a contraction $f : X \to Y$. Let $E \subset X$ be the exceptional set. Let $F$ be any irreducible component of a fiber $f^{-1}(f(x))$ for $x \in E$. Then*

$$\dim F + \dim E \geq \dim X + l(R) - 1.$$

This is not true for contractions of extremal rays on singular varieties.

In particular, Wisniewski's theorem implies that

$$\dim F + \dim E \geq \dim X,$$

which is all we need. I claim that this implies that the contraction $f : X \to Y$ is not IH-small. Indeed, it says that every irreducible component of a fiber of $f : E \to f(E)$ has dimension at least $\dim X - \dim E$, so that $\dim f(E) \leq \dim E - (\dim X - \dim E)$, which translated in terms of codimension says that $\operatorname{codim} E \leq \frac{1}{2} \operatorname{codim} f(E)$. That is, $f$ is not IH-small (Proposition 8.3). □

Now, by Kawamata, if a variety $Y$ has a relative minimal model $f : X \to Y$, then $Y$ also has a relative canonical model $X_0 \to Y$ [14, Th. 3-3-1]. The map $f$ factors $X \to X_0 \to Y$, and we have $K_X = f^*(K_{X_0})$ in $\operatorname{Pic} X \otimes \mathbf{Q}$ by the construction of $X_0$ ($K_X$ is basepoint-free, locally over $Y$).

Thus, if $f : X \to Y$ is an IH-small resolution such that $X$ is projective over $Y$, then $X$ is a relative minimal model of $Y$ by Proposition 8.3, so $Y$ also has a relative canonical model $X_0$, and we have a factorization $X \to X_0 \to Y$ with $K_X = f^* K_{X_0}$. This proves most of Theorem 8.2.

In Theorem 8.2, $Y$ is compact, so $X$ and $X_0$ are also compact. Clearly the map $g : X \to X_0$ is an IH-small resolution and so the last sentence of Section 7 implies that

$$g_*(\chi_i^{n-k}(X)) = \chi_i^{n-k}(X_0) \in H_{2k}(X_0, \mathbf{Q})$$

for all $0 \le k \le n$, $0 \le i \le n-k$. Since $K_X = f^*K_{X_0} \in \text{Pic } X \otimes \mathbf{Q}$, we also have $c_1(X) = f^*c_1(X_0) \in H^2(X, \mathbf{Q})$, and it follows that

$$c_1^k \chi_i^{n-k}(X) = c_1^k \chi_i^{n-k}(X_0) \in \mathbf{Q}$$

(Theorem 8.2). □

## 9. The twisted $\chi_y$ genus

By Section 8, the Chern numbers $c_1^k \chi_i^{n-k}$ can be defined for singular varieties. These Chern numbers can easily be combined into a genus which Höhn called the twisted $\chi_y$ genus: its image is a quotient ring of the elliptic genus quotient ring $\mathbf{Q}[x_1, x_2, x_3, x_4]$ of the bordism ring $MU_* \otimes \mathbf{Q}$. The last paragraph of Höhn's thesis [13] identifies the quotient ring of $MU_* \otimes \mathbf{Q}$ corresponding to the twisted $\chi_y$ genus as

$$\mathbf{Q}[x_1, x_2, x_3, x_4]/(\Delta(x_2, x_3, x_4)).$$

Here we think of the ring $\mathbf{Q}[x_2, x_3, x_4]$ as the ring of Jacobi forms, and $\Delta$ as the discriminant modular form in this ring, which has degree 12 (so that $\Delta \in \mathbf{Q}[x_2, x_3, x_4]$ is the elliptic genus of some linear combination of $SU$-manifolds of complex dimension 12). The proof of this identification of the quotient ring is easy but not quite explicit in Höhn's thesis, so we prove it in this section. The interest of this result from the point of view of this paper is that Section 8 gives an element of this quotient ring associated to any singular variety, when we assume the existence of a relative canonical model.

Let $\chi_y$ be Hirzebruch's $\chi_y$ genus [12], which maps a compact complex $n$-manifold $X$ to the polynomial

$$\chi_y(X) = \chi_0^n(X) + \chi_1^n(X)y + \cdots + \chi_n^n(X)y^n,$$

where we recall that $\chi_i^n(X)$ is the holomorphic Euler characteristic $\chi(X, \Omega^i) = \sum_j (-1)^j \dim H^j(X, \Omega^i)$. It is easy to check that $\chi_y$ is a ring homomorphism $MU_* \to \mathbf{Z}[y]$. To get a homomorphism of graded rings, we redefine $\chi_y(X)$ to be $t^n$ times the above polynomial in $y$ for $X$ of dimension $n$; then $\chi_y$ becomes a homomorphism $MU_* \to \mathbf{Z}[t, y]$ of graded rings, with $t$ in degree 2 and $y$ in degree 0.

As in Section 3, we write $\Lambda_y(E) = \sum_i y^i \Lambda^i E$ for a vector bundle $E$; then the $\chi_y$ genus in $\mathbf{Z}[t, y]$ of a compact complex $n$-manifold can be written as the holomorphic Euler characteristic

$$\chi_y(X) = t^n \chi(X, \Lambda_y(T^*X)).$$

We define the twisted $\chi_y$ genus in $\mathbf{Q}[t, y, z]$ as

$$\chi_{yz}(X) = t^n \chi(X, K_X^{-z} \otimes \Lambda_y(T^*X)).$$

Here $z$ need not be an integer, but the expression still makes sense as a polynomial with rational coefficients by the Hirzebruch-Riemann-Roch theorem. Knowing the twisted $\chi_y$ genus of a manifold is equivalent to knowing all the Chern numbers considered in Section 8.

The twisted $\chi_y$ genus $MU_* \otimes \mathbf{Q} \to \mathbf{Q}[t, y, z]$ is a homomorphism of graded rings with $t$ in degree 1 (corresponding to manifolds of complex dimension 1) and $y$ and $z$ in degree 0. The homomorphism $MU_* \otimes \mathbf{Q} \to \mathbf{Q}[t, y, z]$ is not surjective, but it factors through the complex elliptic genus

$$MU_* \otimes \mathbf{Q} \to \mathbf{Q}[x_1, x_2, x_3, x_4],$$

which is surjective. So we get a well-defined homomorphism:

$$\chi_{yz} : \mathbf{Q}[x_1, x_2, x_3, x_4] \to \mathbf{Q}[t, y, z].$$

We will show that the kernel of this homomorphism is the ideal generated by the discriminant cusp form $\Delta$ in the ring $\mathbf{Q}[x_2, x_3, x_4]$ of Jacobi forms. Thus we can view the twisted $\chi_y$ genus as a surjective homomorphism

$$MU_* \otimes \mathbf{Q} \to \mathbf{Q}[x_1, x_2, x_3, x_4]/(\Delta(x_2, x_3, x_4)),$$

as promised.

To begin with, let $X$ be an $SU$-manifold of complex dimension $n$. Then the elliptic genus $\varphi(X)$ is defined as the power series (see Section 3)

$$\varphi(X) = \chi(X, \prod_{m \geq 1} (\Lambda_{-y^{-1}q^m} T \otimes \Lambda_{-yq^{m-1}} T^* \otimes S_{q^m} T \otimes S_{q^m} T^*)).$$

Here $T$ denotes the virtual bundle $TX - n$ of rank 0. If we define

$$\alpha(X) = \chi(X, \prod_{m \geq 1} (\Lambda_{-y^{-1}q^m} TX \otimes \Lambda_{-yq^{m-1}} T^*X \otimes S_{q^m} TX \otimes S_{q^m} T^*X)),$$

this "unscaled elliptic genus" is related to the usual one by

$$\alpha(X) = \Phi(q, y^{-1})^n \varphi(X),$$

where $\Phi(q, y)$ is the normalization of the Weierstrass sigma function defined in Section 3. The point of introducing the unscaled elliptic genus $\alpha(X)(q, y)$ is that it is evidently related to the $\chi_y$ genus by

$$\begin{aligned} \chi_y(X) &= \alpha(X)(0, -y) \\ &= (1 + y)^n \cdot \varphi(X)(0, -y). \end{aligned}$$

Thus, from Höhn's calculation that there are $SU$-manifolds of complex dimension 2, 3, 4 with elliptic genera $x_2 = 24\mathfrak{p}$, $x_3 = \mathfrak{p}'$, and $x_4 = 6\mathfrak{p}^2 - g_2/2$, we read off from the power series expansions of these Jacobi forms in Section 3 that

$$\begin{aligned} \chi_y(x_2) &= t^2(2 - 20y + 2y^2), \\ \chi_y(x_3) &= t^3(-y + y^2), \\ \chi_y(x_4) &= t^4(-y + 4y^2 - y^3). \end{aligned}$$

Thus the kernel of the $\chi_y$ homomorphism

$$\chi_y : \mathbf{Q}[x_2, x_3, x_4] \to \mathbf{Q}[t, y]$$

is the ideal of "cusp forms" in the ring of Jacobi forms, those which vanish for $q = 0$. If we think of Jacobi forms as sections of powers of the line bundle $\psi_1$ over the universal elliptic curve $\overline{M}_{1,2}$, as in Section 3, a Jacobi cusp form is one which vanishes on the curve $D_0$ (the nodal cubic) in $\overline{M}_{1,2}$, the fiber of $\pi : \overline{M}_{1,2} \to \overline{M}_{1,1}$ over the "cusp" in the modular curve $\overline{M}_{1,1}$. Since the discriminant cusp form $\Delta$ is a section of $12\psi_1$ over $\overline{M}_{1,1}$ with a single zero at the cusp and no other zeros, it pulls back to a section over $\overline{M}_{1,2}$ of $12\psi_1 = 12(\pi^*\psi_1 + D_{0,1})$ whose divisor of zeros is exactly $D_0 + 12D_{0,1}$.

As a result, if we divide a Jacobi form which vanishes on the nodal cubic $D_0$ by $\Delta$, we get a meromorphic Jacobi form which is holomorphic outside $D_{0,1}$. By Lemma 3.1, such a meromorphic Jacobi form is actually holomorphic. Thus every Jacobi form which vanishes on $D_0$ is a multiple of $\Delta$. Equivalently, the kernel of the $\chi_y$ homomorphism

$$\chi_y : \mathbf{Q}[x_2, x_3, x_4] \to \mathbf{Q}[t, y]$$

is the ideal generated by $\Delta$. Yet another way to put this is that the image of the $\chi_y$ genus

$$\chi_y : MSU_* \otimes \mathbf{Q} \to \mathbf{Q}[t, y]$$

is isomorphic to $\mathbf{Q}[x_2, x_3, x_4]/(\Delta)$, where the $SU$-manifolds $x_2, x_3, x_4$ have $\chi_y$ genus as computed above.

It is now easy to determine the image of the $\chi_{yz}$ genus on $MU_* \otimes \mathbf{Q}$. The point is that $MU_* \otimes \mathbf{Q}$ and $MSU_* \otimes \mathbf{Q}$ are both polynomial rings,

$$\begin{aligned} MU_* \otimes \mathbf{Q} &= \mathbf{Q}[x_1, x_2, x_3, \ldots], \\ MSU_* \otimes \mathbf{Q} &= \mathbf{Q}[x_2, x_3, x_4, \ldots], \end{aligned}$$

and a generator for $MSU_* \otimes \mathbf{Q}$ in any complex dimension $n \geq 2$ is also a generator for $MU_* \otimes \mathbf{Q}$ in the same dimension. (Novikov [29] proved this as well as more precise integral information.) Thus we can say that $MU_* \otimes \mathbf{Q}$ is generated as a $\mathbf{Q}$-algebra by $\mathbf{CP}^1$ together with the image of $MSU_*$.

For a complex manifold with trivial canonical bundle, the $\chi_{yz}$ genus is equal to the $\chi_y$ genus (that is, it is a polynomial only in $y$, not involving $z$). So, by the previous section, we know that the image of the $\chi_{yz}$ genus on the image of $MSU_* \otimes \mathbf{Q}$ in $MU_* \otimes \mathbf{Q}$ is the ring

$$\mathbf{Q}[x_2, x_3, x_4]/(\Delta(x_2, x_3, x_4)).$$

On the other hand, the $\chi_{yz}$ genus of $\mathbf{CP}^1$ is $1 - y + 2z$. Since this involves $z$, it does not satisfy any relations with $x_2$, $x_3$, and $x_4$. Thus the image of the $\chi_{yz}$ genus on $MU_* \otimes \mathbf{Q}$ is the ring

$$\mathbf{Q}[x_1, x_2, x_3, x_4]/(\Delta(x_2, x_3, x_4)),$$

where $x_1 = t(1-y+2z)$ and $x_2, x_3, x_4$ are as above. Explicitly, we can expand the discriminant cusp form $\Delta$ in terms of the Jacobi forms $x_2 = 24\mathfrak{p}$, $x_3 = \mathfrak{p}'$, and $x_4 = 6\mathfrak{p}^2 - g_2/2$ defined in Section 3 by

$$\begin{aligned}\Delta &= g_2^3 - 27g_3^2 \\ &= g_2^3 - 27(4\mathfrak{p}^3 - g_2\mathfrak{p} - (\mathfrak{p}')^2)^2 \\ &= -\frac{1}{32}x_2^3x_3^2 + \frac{1}{16}x_2^2x_4^2 + \frac{9}{2}x_2x_3^2x_4 - 27x_3^4 - 8x_4^2.\end{aligned}$$

Theorems 8.1 and 8.2 assert the possibility of defining an element of the above ring associated to any singular variety, when we assume the existence of the relative canonical model.

There is a different way to describe this calculation of the image of the $\chi_y$ genus on $SU$-manifolds, rationally: the only linear relations satisfied by the $\chi_y$ genus of an $SU$-manifold are those coming from Serre duality together with the one other relation found by Libgober and Wood [22]. From Serre duality, the $\chi_y$ genus of an $SU$-manifold of complex dimension $n$ is a polynomial $\chi(y)$ of degree $n$ such that $\chi(1/y) = (-1/y)^n\chi(y)$. Also by Serre duality, if $n$ is odd, the Todd genus $\chi(0)$ is 0. Libgober and Wood's relation, which involves the Euler characteristic $\chi(-1)$, says that

$$\chi''(-1) - \frac{n(3n-5)}{12}\chi(-1) = 0.$$

The fact that there are no other rational linear relations satisfied by the $\chi_y$ genus of an $SU$-manifold is proved by computing the dimension of the ring $\chi_y(MSU_*) \otimes \mathbf{Q} = \mathbf{Q}[x_2, x_3, x_4]/(\Delta)$ in each degree.

Cambridge University, Cambridge CB2 1SB, England
*E-mail address*: b.totaro@dpmms.cam.ac.uk

## References

[1] M. Atiyah, On analytic surfaces with double points, *Proc. Roy. Soc. London*, Ser. A **247** (1958), 237–244.

[2] N. Baas, On bordism theory of manifolds with singularities, *Math. Scand.* **33** (1973), 279–302.

[3] G. Barthel, J.-P. Brasselet, K.-H. Fieseler, O. Gabber, and L. Kaup, Relèvement de cycles algébriques et homomorphismes associés en homologie d'intersection, *Ann. of Math.* **141** (1995), 147–179.

[4] P. Baum, W. Fulton, and R. MacPherson, Riemann-Roch and topological $K$-theory for singular varieties, *Acta Math.* **143** (1979), 155–192.

[5] J.-P. Brasselet and G. Gonzalez-Springberg, Sur l'homologie d'intersection et les classes de Chern des variétés singulières (espaces de Thom, exemples de J.-L. Verdier et M. Goresky), with an appendix by Verdier, in *Travaux en Cours*, 23. *Geometrie Algébrique et Applications*, II (*La Rábida*, 1984), 5–14, Hermann, Paris (1987).

[6] P. Deligne, Courbes elliptiques: formulaire (d'après J. Tate), *Modular Functions of One Variable* IV, *Lecture Notes in Math.* **476**, 53–73, Springer-Verlag, Berlin (1975).

[7] P. du Bois, Complexe de de Rham filtré d'une variété singulière, *Bull. Soc. Math. France* **109** (1981), 41–81.
[8] M. Eichler and D. Zagier, *The Theory of Jacobi Forms*, Birkhäuser, Boston (1985).
[9] W. Fulton, *Intersection Theory*, Springer-Verlag, New York (1984).
[10] M. Goresky and R. MacPherson, Intersection homology, II, *Invent. Math.* **72** (1983), 77–130.
[11] ———, Problems and bibliography on intersection homology, in *Intersection Cohomology*, ed. A. Borel et al., 221–233, Birkhäuser, Boston (1984).
[12] F. Hirzebruch, *Topological Methods in Algebraic Geometry*, 3rd ed., Springer-Verlag, New York (1966).
[13] G. Höhn, Komplexe elliptische Geschlechter und $S^1$-äquivariante Kobordismustheorie, Diplomarbeit, Bonn, Germany (1991). Also at: http://baby.mathematik.uni-freiburg.de/papers/ .
[14] Y. Kawamata, K. Matsuda, and K. Matsuki, Introduction to the minimal model program, in *Algebraic Geometry* (*Sendai*, 1985), ed. T. Oda, *Adv. Stud. Pure Math.* **10**, 283–360, Kinokuniya–North Holland (1987).
[15] F. Knudsen, Projectivity of the moduli space of stable curves, II, *Math. Scand.* **52** (1983), 1225–1265.
[16] J. Kollár, The structure of algebraic threefolds: an introduction to Mori's program, *Bull. AMS* **17** (1987), 211–273.
[17] ———, Flips, flops, minimal models, etc., in *Surveys in Differential Geometry* (*Cambridge*, *MA*, 1990), 113–199, Lehigh Univ., Bethlehem, Pa. (1991).
[18] M. Kreck and S. Stolz, $\mathbf{HP}^2$-bundles and elliptic homology, *Acta Math.* **171** (1993), 231–261.
[19] I. Krichever, Generalized elliptic genera and Baker-Akhiezer functions, *Math. Notes* **47** (1990), 132–142.
[20] P. Landweber, ed., *Elliptic Curves and Modular Forms in Algebraic Topology*, *Lecture Notes in Math.* **1326**, Springer-Verlag, New York (1988).
[21] S. Lang, *Elliptic Functions*, Springer-Verlag, New York (1987).
[22] A. Libgober and J. Wood, Uniqueness of the complex structure on Kähler manifolds of certain homotopy types, *J. Differential Geom.* **32** (1990), 139–154.
[23] J. Milnor, On the cobordism ring $\Omega^*$ and a complex analogue, I, *Amer. J. Math.* **82** (1960), 505–521.
[24] J. Morava, A product for the odd-primary bordism of manifolds with singularities, *Topology* **18** (1979), 177–186.
[25] S. Mori, Threefolds whose canonical bundles are not numerically effective, *Ann. of Math.* **116** (1982), 133–176.
[26] ———, Flip theorem and the existence of minimal models for threefolds, *J. AMS* **1** (1988), 117–253.
[27] D. Mumford, Towards an enumerative geometry of the moduli space of curves, in *Arithmetic and Geometry*, ed. Michael Artin and John Tate, 271–328, Birkhäuser, Boston (1983).
[28] M. Narasimhan and S. Ramanan, Generalized Prym varieties as fixed points, *J. Ind. Math. Soc.* **39** (1975), 1–19.
[29] S. P. Novikov, Homotopy properties of Thom complexes (Russian), *Mat. Sb.* **57** (1962), 407–442.
[30] I. Porteous, Blowing up Chern classes, *Proc. Camb. Phil. Soc.* **56** (1960), 118–124.
[31] M. Reid, Decomposition of toric morphisms, in *Arithmetic and Geometry*, v. II, 395–418. Birkhäuser, Boston (1983).
[32] M. Saito, Modules de Hodge polarisables, *Publ. Res. Inst. Math. Sci.* **24** (1988), 849–995.
[33] ———, Mixed Hodge modules, *Publ. Res. Inst. Math. Sci.* **26** (1990), 221–333.
[34] J. Wisniewski, On contractions of extremal rays of Fano manifolds, *J. Reine Angew. Math.* **417** (1991), 141–157.